\documentclass[12pt]{amsart}
\usepackage{amssymb}
\numberwithin{equation}{section}
\theoremstyle{plain}
\newtheorem{theorem}{Theorem}[section]
\newtheorem{prop}[theorem]{Proposition}
\newtheorem{corollary}[theorem]{Corollary}
\newtheorem{lemma}[theorem]{Lemma}

\theoremstyle{definition}

\theoremstyle{remark}

\newcommand{\ba}{\begin{eqnarray}}
\newcommand{\be}{\begin{equation}}
\newcommand{\ea}{\end{eqnarray}}
\newcommand{\ee}{\end{equation}}
\newcommand{\benn}{\begin{equation*}}
\newcommand{\eenn}{\end{equation*}}

\def\R{{\mathbb R}}
\def\N{{\mathbb N}}
\def\Z{{\mathbb Z}}
\def\Cpx{{\mathbb C}}

\def\supp{\operatorname{supp}}
\def\tr{\operatorname{Tr}}
\def\meas{\operatorname{meas}}

\def\Rone{{\mathbb R}}
\def\Rn{{{\mathbb R}^n}}

\def\Fcal{\mathcal{F}}
\def\Dcal{\mathcal{D}}
\def\Lcal{\mathcal{L}}

\def\diag{{\;{\rm diag}\;}}
\def\jp#1{{\left\langle{#1}\right\rangle}}

\title{Regularity properties, representation of solutions and
spectral asymptotics of systems with multiplicities}

\author[]{Ilia Kamotski and Michael Ruzhansky}
\address{
  Ilia Kamotski,
  Michael Ruzhansky:
  \endgraf
  Department of Mathematics
  \endgraf
  Imperial College of Science, Technology and Medicine
  \endgraf
  180 Queen's Gate, London SW7 2BZ, UK
  \endgraf
  {\it E-mail addresses} {\rm kamotski@ic.ac.uk, ruzh@ic.ac.uk}
  \endgraf
  \medskip
  }
\date{\today}
\begin{document}
\maketitle

\begin{abstract}
Properties of solutions of generic hyperbolic systems with multiple
characteristics with diagonalizable
principal part are investigated. Solutions are represented as a Picard
series with terms in the form of 
iterated Fourier integral operators. It is shown that this series is
an asymptotic expansion with respect to smoothness under quite
general geometric conditions.
Propagation of singularities and sharp regularity properties
of solutions are obtained. Results are applied to establish regularity
estimates for scalar weakly hyperbolic equations with involutive 
characteristics. They are also applied to derive the first and
second terms of spectral
asymptotics for the corresponding elliptic systems.  
\end{abstract}

\section{Introduction}

Let $X$ be a smooth manifold without boundary of dimension $n\geq 3.$
Let $P$ be an elliptic self-adjoint pseudo-differential operator
of order one acting on half-densities on $m$-dimensional
crossections of vector
bundles on $X$.
We consider the following Cauchy problem for $u=u(t,x)$

\begin{equation}
    \left\{
      \begin{array}{ll}
         iu^\prime-Pu=0, & (t,x)\in [0,T]\times X,
          \\
         u|_{t=0}=u_0.
      \end{array} \right.
      \label{eq:eq1}
\end{equation}
It is well known that if equation (\ref{eq:eq1}) is strictly hyperbolic,
the system can be diagonalized and its solution can be given as a sum of
Fourier integral operators applied to Cauchy data (e.g. \cite{Duis}).
An important question
that has been studied over many years is what happens when $P$ has
multiple characteristics.

Since we will be mostly interested in local properties of
solutions, we may already assume that $P$ acts
on functions, and can think of it as an $m\times m$ matrix of pseudo-differential
operators of order one and we think of $u_0$ as of an $m$-vector.

Let $A(x,\xi)$ denote the principal symbol of $P$. If $A$ is a diagonal matrix,
properties of system (\ref{eq:eq1}) have been studied by many authors.
For example, in \cite{KTT} and \cite{Kumano-go-book} Kumano-go and
coauthors used the calculus of Fourier
integral operators with multi-phases to show that the Cauchy problem
(\ref{eq:eq1}) is well-posed in $L^2$, Sobolev spaces $H^s$, and to study
its propagation of singularities. Systems with
symmetric principal part $A$ have been extensively studied as well (e.g.
\cite{Kucherenko}, \cite{Ivriibook}, etc.) In a generic situation, such systems
have double characteristics, and their normal forms have been found by
Braam and Duistermaat \cite{BD}.
Recently, Colin de Verdiere \cite{ColindV} used
these representations to derive some asymptotic properties of such systems.
Polarization properties of similar systems were studied by
Dencker in \cite{Dencker}.

More elaborate analysis of system (\ref{eq:eq1}) becomes possible if
we assume that the principal symbol matrix $A(x,\xi)$ is smoothly
microlocally diagonalizable
with smooth eigenvalues $a_j(x,\xi)$ and smooth eigenspaces.
Then, as it was pointed
out by Rozenblum in \cite{Rozenblum}, there exists a finite dimensional cover
$\tilde{X}$ of $X$
such that $A$ lifted to $\tilde{X}$ can be
globally digonalized provided $X$ is compact.
In this situation
Rozenblum showed that the Picard series for this problem gives an expansion
with respect to smoothness in the case of non-involutive characteristics.
In other words, one assumes that
if $a_j(x,\xi)=a_k(x,\xi)$ for $j\not=k$, then the Poisson bracket
$\{a_j,a_k\}(x,\xi)\not=0.$
This means that at all points of multiplicity,
bicharacteristics intersect characteristics surfaces transversally.
However, this condition is non-generic even for
diagonal systems. For example, it is clear that if one of characteristics
has a fold, there may be a point where this transversality condition
fails, and it is not possible to remove it by small perturbations.

One purpose of this paper is to present results removing
the transversal intersection condition.
We will allow characteristics to be involutive of finite type and some
characteristics to be identically equal.
Operators satisfying our
Condition C below will be also generic in the class of smoothly microlocally
diagonalizable systems. Below we will explain that the microlocal diagonalizability
condition is quite natural when considering weakly hyperbolic scalar
equations with Levi conditions (Examples 1 and 2).
This is also the case for Maxwell equations (e.g. \cite{BS}).

We will
investigate regularity properties of system (\ref{eq:eq1}) in this generic
setting.
Even in the case when the system is strictly hyperbolic, $L^p$
properties
of solutions have been studied for many years, since already this case has
several important applications for nonlinear equations and harmonic analysis.
Regularity properties of non-degenerate Fourier integral operators have been
established by Seeger, Sogge and Stein in \cite{SSS}. They showed that
a Fourier integral operator $T$
of order zero satisfying local graph condition, is locally bounded from
$(L^p_\alpha)_{comp}$ to $L^p_{loc}$
for $1<p<\infty$ and $\alpha=(n-1)|1/p-1/2|$.
As a consequence they showed that if system \eqref{eq:eq1} is strictly
hyperbolic, there is a loss of $\alpha$ derivatives in $L^p$, i.e.
$u_0\in L^p_\alpha$ implies $u(t,\cdot)\in L^p$. Moreover, if at least
one of characteristic
roots $a_j$ is elliptic, the loss of $\alpha$ derivatives is sharp.
If none of $a_j$'s is elliptic, this result can be improved (\cite{RuzhSurv}).

Our Theorem \ref{th:lp} will establish a similar property for systems
\eqref{eq:eq1} with multiplicities. Moreover, this will imply $L^p$
estimates for scalar weakly hyperbolic equations with involutive
characteristics. It is known that in general weakly hyperbolic
cases one often loses regularity even in $L^2$.  However,
in the case of involutive characteristics
the equation can be diagonalized and in Theorem \ref{th:inv}
we will give $L^p$ estimates for such equations.
This will, on one hand, extend the $L^p$ result of
Seeger, Sogge and Stein to weakly hyperbolic equations and systems
with multiplicities, while on the other hand establishing $L^p$ estimates
for systems considered by Kumano-go, Rozenblum, and others.
The result will be general for scalar weakly hyperbolic equations satisfying
Levi conditions with characteristics satisfying Condition C below. 

Note that if a scalar operator strictly hyperbolic and we write it in the
form \eqref{eq:eq1}, we can diagonalize $P$ together with lower order
terms (e.g. \cite{Kumano-go-book}) and split it into $m$ scalar equations,
for which many things are known. However, in the case of multiple
characteristics this is impossible, so a more elaborate analysis is needed.

Now we will formulate our main assumption.
Let us define operator
$$ H_{a_j} f=\{a_j,f\}, j=1,\ldots, m,$$
where
$$H_g(f)=\{g,f\}=\sum_{k=1}^n \left(\frac{\partial g}{\partial \xi_k}
\frac{\partial f}{\partial x_k}-\frac{\partial g}{\partial x_k}
\frac{\partial f}{\partial \xi_k}\right)$$
is the usual Poisson bracket.
Our assumption is that at points of multiplicity $a_j=a_k$ of non-identical
characteristics $a_j$ and $a_k$,
bicharacteristics
of $a_j$ intersect level sets $\{a_k=1\}$ with finite order, i.e.
$H_{a_j}^M a_k\not=0$ for some $M$, at points where $a_j=a_k$.
In other words, we allow involutive characteristics of finite type,
and formulate our main condition.

\bigskip
\noindent
Condition C: \\
Suppose that there exists $M\in\N$ such that
if for some $j$ and $k$, $a_j$ and $a_k$ are not identically the same, then
\begin{equation}
\begin{array}{l} a_j(x,\xi)=a_k(x,\xi), \ (x,\xi)\in T^*X 
\Longrightarrow \\
H_{a_j}^\lambda
a_k(x,\xi)=\{a_j,\{a_j,\cdots\{a_j,a_k\stackrel{\lambda}{\overbrace
{\}\}\ldots\}}}(x,\xi)\not=0, \end{array}
\label{eq:eq3}
\end{equation}
for some $\lambda\leq M$. While the function $M=M(x,\xi)$ is locally bounded,
it is allowed to grow at infinity.

We note here that the transversality assumption
of Rozenblum \cite{Rozenblum} requires (\ref{eq:eq3}) to hold with $M=1$.
Strictly hyperbolic case is also covered by this condition (in which
case we set $M=0$). The case of $a_j$ and $a_k$ defining glancing hypersurfaces
(as in Melrose \cite{Melrose}) corresponds to $M=2$.

In Section 2 we will give several examples of characteristics
satisfying condition C, in particular those arising from weakly hyperbolic
scalar equations with involutive characteristics.
Such equations and propagation of their singularities have been
analyzed in
\cite{Chazarain}, \cite{MO}, \cite{Morimoto}, \cite{Ichinose}, \cite{IM},
etc.
We will also establish estimates in $L^p$ and other spaces for
the weakly hyperbolic equations or systems satisfying condition C.

It is interesting to note that conditions similar to Condition C appeared
in the study of subelliptic operators (e.g. H\"ormander \cite[Chapter 27]{Hor}).
For instance, in the case of $2\times 2$ systems $P$ with characteristics
$a_1$ and $a_2$, we can consider operators $Q$ with principal symbol
$q=a_1+ia_2$. Then microlocal subellipticity of $Q$ implies Condition C,
with some $M$ dependent on the loss of regularity for $Q$, which, therefore,
implies the Weyl formula for $P$ (Theorem \ref{th:spectral}), regularity
estimates for \eqref{eq:eq1} and all other results of this paper.

Now we will give an informal explanation of the strategy of our analysis.
First, let us follow \cite{Rozenblum} to show that microlocal
diagonalizability implies a local one
on some cover $\tilde{X}$ of $X$ with finitely many leaves.
For this argument we assume that $X$ is compact. Since all the 
analysis of this
paper will be local, if $X$ is not compact, we can always assume that
the amplitude of $P(x,D)$ is compactly supported.

Let an elliptic pseudo-differential operator $P(x,D)$ of order one
act on sections of an $m$-dimensional Hermitian vector bundle $E$.
Let $L^2(E)$ be the space of sections of half-densities on $E$
and let $P$ be self-adjoint on $L^2(E)$. Let $E^\prime$ be the lifting
of $E$ to $T^*X$. Then for each $(x,\xi)\in T^*X$ the principal symbol
$A(x,\xi)$ of $P(x,D)$ is a Hermitian isomorphism of sections of
$E^\prime$. Without loss of generality we can assume that $A(x,\xi)$
is positive definite. Indeed, if it has both positive and negative
eigenvalues, it is possible to globally block-diagonalize $A(x,D)$
with some suitably chosen pseudo-differential operators, to reduce
it to a direct sum of positive and negative definite operators.
Then each of these operators can be analyzed independently.

Let us assume that the principal symbol $A(x,\xi)$ is microlocally
diagonalizable. This means that microlocally in
$\Lambda\subset T^*X$ such that
$E^\prime|_\Lambda\cong \Lambda\times\Cpx^m$, principal symbol
$A(x,\xi)|_\Lambda$ has $m$ smooth non-negative eigenvalues
$a_j(x,\xi)$ and one dimensional eigenspaces $V_j(x,\xi)$, and such
diagonalizations are compatible in intersecting cones.
In this situation Lemma \ref{l:diag} insures that there is a
cover $\tilde{X}$ of $X$ with finitely many leaves such that
the principal symbol of the lifting of $P(x,D)$ to $T^*\tilde{X}$
can be globally diagonalized modulo lower order terms.
Note that since dimensions of $X$ and $\tilde{X}$ are the same and
because of formula \eqref{eq:cover} all our results on $\tilde{X}$
will imply corresponding results on $X$.
Therefore, we may assume that the principal symbol $A$ of operator $P$
may be smoothly diagonalized over compact subsets of $X$,
that is
$$ P=A+B, A={\rm diag}\{A_1,\ldots,A_m\},$$
where
$A_j\in \Psi^1$ are scalar pseudo-differential operators
 with principal symbols $a_j(x,\xi)$
(eigenvalues of $A$). Here
$a_j$'s may be identically equal to each other
or may intersect with any finite order, according to our
Condition C.
Here $B$ is an $m\times m$ matrix of
pseudo-differential operators or order zero. We may also assume
that $B_{jj}=0$ for $1\leq j\leq m$, if we add these terms to the
diagonal of $A$.
For the moment we will also assume that none of $a_j$'s are
identical.
Otherwise, if some of $a_j$'s
being identically equal to each other locally at some points,
the construction is slightly different, but all the results remain valid.
This will be carried out in detail in Section 3 in the proof of Theorem
\ref{th:lp}.
Substitution $U=e^{-iAt}V$ leads to the equation
 \begin{equation}
       \left\{ \begin{array}{l}
         V^\prime= Z(t) V,
          \\
         V|_{t=0}=I,
      \end{array} \right.
      \label{eq:eq4}
\end{equation}
with $Z(t)=-ie^{iAt} B e^{-iAt}.$
Writing the Picard series for problem (\ref{eq:eq4}), we obtain the expansion
\begin{equation}
  \begin{array}{lll}
 V(t)&=&I+\int_0^t Z(t_1) dt_1 + \int_0^t\int_0^{t_1} Z(t_1) Z(t_2) dt_2
dt_1+\cdots
 \end{array}
 \label{eq:series}
\end{equation}
A general term of this series is
\be
Q_l=\stackrel{l}{\overbrace{\int_0^t\int_0^{t_1}\cdots
\int_0^{t_{l-1}}}} e^{iA_{j_1}t_1}b_{j_1 j_2}e^{iA_{j_2}(t_2-t_1)}\cdots
dt_l\ldots dt_1.
\label{eq:ql}
\ee
It is easy to see that $||Q_l||_{L^2\to L^2}\leq C/l!$
and that series (\ref{eq:series}) converges in $L^2$ and in $H^s$.
Using the notion of a multi-phase for Fourier integral operators,
Kumano-go et al. (\cite{KTT}, \cite{Kumano-go-book})
studied propagation of singularities of $Q_l$.
Instead of introducing multi-phases for
$Q_l$, we will analyze operators $Q_l$ in more detail and will
show its smoothing properties in Sobolev spaces under Condition C.
Here, contrary to the transversal non-involutive case of
Rozenblum (when $M=1$),
we do not have good control on the singular supports of
integral kernels of operators $Q_l$, so more elaborate geometric
analysis is required.

In fact, our Theorem \ref{th:q} asserts
that $Q_l(t)$ maps $L^2$ to some Sobolev space
$H^{p(l)}$, and that $p(l)\to\infty$
as $l\to\infty.$ Then we will show that this allows to treat the Picard
expansion as a series with finitely many terms, with many conclusions,
such as $L^p$ estimates for solutions and spectral asymptotics of $P$.

In particular, $L^p$ estimates will follow from a general principle which
we will prove in Theorem \ref{th:regularity} for equation
$u^\prime-Z(t)u=f, u(0)=u_0$. In Corollaries \ref{cor:c0} and \ref{cor:c1}
we will show regularity of solutions of this equation for
pseudo-differential operators of order zero $Z(t)\in\Psi^0$ or 
for Fourier integral operators of negative orders $Z(t)\in I^{-\epsilon}$.
In general, there may be problems with this construction for $Z(t)\in I^0$,
but there we can use the structural properties of $Z(t)$ in 
\eqref{eq:eq4}. 

Everywhere in this paper $\Psi^\mu=\Psi^\mu_{1,0}(X)$ will denote the
space of classical pseudo-differential operators of order $\mu$ of type $(1,0)$.
The space of Fourier integral operators of order $\mu$ with amplitudes
of type $(1,0)$ will be denoted by $I^\mu$. All Fourier integral operators
in this paper will be non-degenerate, which means that its canonical
relation satisfies the local graph condition, i.e. it is a graph of a
symplectic diffeomorphism from $T^*X$ to itself. Constants $C$ may be
different in different formulas throughout this paper. We will use
the following notation for norms and spaces. By $L^p_\alpha$ we will
denote the Sobolev space of functions $f$ such that
$(I-\Delta)^{\alpha/2}f\in L^p$.
For a function $f$ we will denote its $L^p$-norm by $\|f\|_{L^p}$
and its
Sobolev $H^s=L^2_s$ norm by $\|f\|_s$. If $T$ is an operator, by
$\|T\|_{s}$ we will denote its operator norm from $L^2$ to $H^s$.

Results of this paper can be established also in the case of
operators $P$ dependent on $t$. This will be the subject of
a separate paper since it will also involve the analysis of 
non-smooth coefficients.

We would like to thank Ari Laptev for drawing our attention to
Rozenblum's paper \cite{Rozenblum}, Grigori Rozenblum for several
valuable references and Yuri Safarov for discussions. The research
was supported by EPSRC grant GR/R67583/01.

\section{Main results}

Now we will formulate our results concerning terms of
the Picard series \eqref{eq:series}
and solutions to systems (\ref{eq:eq1}) and (\ref{eq:eq4}).
We will also give a Weyl formula for $P$.

Our first main result will be Theorem \ref{th:q} on the smoothing
properties of terms $Q_l$. The other two important results will be
Theorem \ref{th:lp} on the $L^p$-regularity of solutions to system
(\ref{eq:eq1}) and Theorem \ref{th:spectral} on the spectral asymptotics
for elliptic operator $P$ satisfying Condition C.
To obtain $L^p$-estimates, we use Theorem \ref{th:regularity},
which we regard as a general principle behind regularity estimates
for general Cauchy problems based on several natural properties of
the right hand side operators $Z(t)$. We will illustrate its use
in several situations in Corollaries 1-3. Theorem
\ref{th:singularities} is a statement on the propagation
of singularities of operators $Q_l$ or, more generally, of
Fourier integral operators in which the frequency integration is
performed over a cone rather than over the whole space.
Note that
these results will hold for microlocally smoothly block-diagonalizable
operators with any (finite) geometry of characteristics, i.e. characteristics
satisfying our Condition C.

The following theorem establishes a smoothing property of operators
$Q_l$ under Condition C.

\begin{theorem}
Let condition C be satisfied, that is
assume that there is $M$ such that
for any $(x,\xi)\in T^*X$ and any $1\leq j,k\leq m$, $j\not=k,$ with
$a_j$ and $a_k$ not identically the same and
$a_j(x,\xi)=a_k(x,\xi)$, we have
$$\{a_j,\{a_j,\cdots\{a_j,a_k\stackrel{\lambda}{\overbrace
{\}\}\ldots\}}}(x,\xi)\not=0 $$
for some number $\lambda\leq M$.
Then for sufficiently large $l$ operator $Q_l$ in \eqref{eq:ql} is
bounded from $L^2_{comp}$ to $H^N_{loc}$, where
$$N=\frac{((-3n/2-2)(3[l/2]-1-n)+([l/2]-n-1)(l/(2M)-n-1))}{(3[l/2]-2n-2+l/(2M))}
\sim \frac{l}{6M+2}.$$
\label{th:q}
\end{theorem}
Note that the exact order $N$ can be improved. However, it is most important
that it increases to infinity as $l\to\infty$.
This theorem implies, in particular, that the series \eqref{eq:series},
i.e. the series
$$V(t)=I+Q_1(t)+Q_2(t)+\cdots$$
is a series with respect to smoothness. This fact allows one to refine
the study of propagation of singularities and regularity properties of
solutions to systems (\ref{eq:eq1}) and (\ref{eq:eq4}).

It turns out that smoothing properties of $Q_l$ in Theorem \ref{th:q} can
be used to establish local $L^p$ properties of solutions to systems with
multiplicities (\ref{eq:eq1}). For strictly hyperbolic
equations such estimates have been established by Seeger, Sogge and Stein
in \cite{SSS}, and some optimal estimates were given in \cite{RuzhSurv}.
Our next result concerns regularity of solutions to Cauchy
problem \eqref{eq:eq1}. In this theorem we also allow the lower
order term $B$ to depend on $t$.

\begin{theorem}\label{th:lp}
Let $1<p<\infty$ and $\alpha=(n-1)|1/p-1/2|.$ Let $P=P(t,x,D_x)$
be an $m\times m$ matrix of elliptic classical
pseudo-differential operators of order one.
Let $$P(t,x,D_x)=A(x,D_x)+B(t,x,D_x),$$
where $A$ is a symmetric matrix of pseudo-differential operators of
order one and $B$ is a
matrix of operators of order zero. Assume that the matrix $A$ is
smoothly (microlocally) diagonalizable, with smooth eigenspaces and
real eigenvalues $a_j(x,\xi)$,
satisfying condition C. Then for any compactly
supported $f\in L^p_\alpha\cap L^2$, the solution $u=u(t,x)$ of
the Cauchy problem
\begin{equation}\label{eq:cauchya}
  i\frac{\partial u}{\partial t}-P(t,x,D_x)u = 0, \;\;
  u(0)=f,
\end{equation}
satisfies $u(t,\cdot)\in L^p_{loc}$ for all $0<t\leq T.$
Moreover, there is a constant $C>0$ such that
\[
\sup_{0\leq t\leq T} ||u(t,\cdot)||_{L^p}\leq C_T||f||_{L^p_\alpha}.
\]
\end{theorem}
We note that it is sufficient to only assume that $A$ is microlocally 
diagonalizable. Then by Lemma \ref{l:diag} we can first have the
statement of Theorem \ref{th:lp} on the cover $\tilde{X}$. Then, using
formula \eqref{eq:cover} and the fact that $\dim\tilde{X}=\dim X$, we
get the same conclusion on $X$.

As a consequence, if the Cauchy data $u_0$ is compactly supported,
we obtain local estimates in other spaces as well:
\begin{itemize}
\item $u_0\in L^p_{s+\alpha}$ implies $u(t,\cdot)\in L^p_s$, $s\in\Rone.$
\item $u_0\in {\rm Lip}(s+(n-1)/2)$ implies
$u(t,\cdot)\in {\rm Lip}(s).$
\item Let $1<p\leq q\leq 2.$ Then
$u_0\in L^p_{s-1/q+n/p-(n-1)/2}$ implies $u(t,\cdot)\in L^q_s.$
Dual result holds for $2\leq p\leq q<\infty.$
\end{itemize}
The proof of Theorem \ref{th:lp} is based on the other main result
Theorem \ref{th:regularity} which we will discuss in the next section.
Estimates in other spaces follow by standard methods of harmonic analysis
(\cite{St}).

Now we will give some examples where condition C holds while
the tranversality condition ($M=1$) fails. This is, for example,
the case when pairs $a_j, a_k$ define glancing hypersurfaces or
when we consider Maxwell systems with variable coefficients.
Below we will concentrate on systems arising from scalar weakly
hyperbolic equations with Levi conditions.

\medskip
\noindent
{\bf Example 1.} In scalar equations with Levi conditions studied
by Chazarain \cite{Chazarain}, Mizohata-Ohya \cite{MO},
Zeman \cite{Zeman}, one assumed that
$\{a_j,a_k\}=C_{jk} (a_j-a_k).$ It is clear that in this situation
$a_j(x,\xi)=a_k(x,\xi)$ implies $\{a_j,a_k\}(x,\xi)=0$.
However, in a general case when $C_{jk}(x,\xi)$ is non-constant,
condition C is satisfied generically.

\medskip
\noindent
{\bf Example 2.} Let $L$ be a scalar operator with involutive
characteristics. More precisely, let
us denote $\partial_j=D_t+\lambda_j(t,x,D_x)$ and let
\be
L=\partial_1\cdots\partial_m+\sum_{k<m}b_{j_1,\cdots,j_k}\partial_{j_1}
\cdots\partial_{j_k}+c,
\label{eq:inv1}
\ee
where $b(t,x,D_x), c(t,x,D_x)\in \Psi^0$ are pseudo-differential operators
of order zero for all $t\in [0,T]$.
We will assume that symbols of all operators are infinitely differentiable with
respect to $t$ in the topology of symbols of the corresponding order.
Let us assume that operator $L$ has involutive characteristics, i.e. that
$$
[\partial_j,\partial_k]\equiv \partial_j\partial_k-\partial_k\partial_j=
\alpha_{jk}\partial_j+\beta_{jk}\partial_k+\gamma_{jk},
$$
where $\alpha_{jk},\beta_{jk},\gamma_{jk}\in\Psi^0$ are
pseudo-differential operators of order zero.
Then it was shown by Morimoto in \cite{Morimoto} that
the Cauchy problem for the equation $Lu=f$ is diagonalizable (with
$1+\sum_{j=1}^{m-1} m!/j!$ components).
Even in the simplest case of characteristics not depending on $t$, we have
$$
\{\lambda_j,\lambda_k\}=\alpha_{jk}(x,\xi)(\lambda_j-\lambda_k)+
\gamma_{jk},
$$
similar to Example 1.

Propagation of singularities of systems with vanishing Poission brackets
has been also studied in these situations. For example,
Iwasaki and Morimoto \cite{IM} studied propagation of singularities of
$3\times 3$ systems, where the second Poisson bracket vanish. Also,
Ichinose \cite{Ichinose} studied $2\times 2$ systems with vanishing second
Poisson brackets. Theorem \ref{th:lp} implies a precise statement on
$L^p$ estimates.

\begin{theorem}
Let $1<p<\infty$, $\alpha=(n-1)|1/p-1/2|$ and $s\in\R$.
Let $L$ be as in \eqref{eq:inv1} and suppose that principal symbols
$a_j(x,\xi)$ of $\lambda_j$ satisfy condition C and do not depend on $t$.
Let $u$ be a solution to the Cauchy problem
\be
\left\{
\begin{array}{l}
Lu=0, \\
\partial_t^j u(0,x)=g_j(x), \ 0\leq j\leq m-1,
\end{array}
\right.
\label{eq:invs}
\ee
and let Cauchy data $g_j\in L^p_{\alpha-j+s}$ be compactly supported.
Then 
$u(t,\cdot)\in (L^p_s)_{loc}$ for all $t\in [0,T]$ and
\be
\sup_{t\in [0,T]} \| \partial^J u(t,\cdot)\|_{L^p_s}\leq
C\sum_{j=0}^{m-1} \| g_j\|_{L^p_{\alpha-j+s+m-1}},
\label{eq:inve}
\ee
where $\partial^J=\partial_{j_1}\ldots\partial_{j_k}$, $k\leq m-1$, and 
$(j_1,\ldots,j_k)$ being permutations of some elements of $\{1,\ldots,m\}$.
\label{th:inv}
\end{theorem}
Note that in the strictly hyperbolic case as well as in some very special
cases of operator $L$ in \eqref{eq:inv1} (e.g. when all $b$ and $c$ are 
zero), following the method described by Treves in \cite{TrevesFIO}
and estimates for Fourier integral operators,
it is possible to obtain the estimate for the Sobolev norm 
$\|u\|_{L^p_{s+m-1}}$ in the left hand side of \eqref{eq:inve}. 

Let us now give a final example of $L^p$ estimates for second order equations,
which we will prove in the next section.

\medskip
\noindent
{\bf Example 3.}
Let us consider the second order equation
$$u^{\prime\prime}+b(x,D_x)u^\prime+c(x,D_x)u=0,$$
where $b\in\Psi^1, \; c\in\Psi^2.$
Let us denote $\langle x\rangle=(1+x^2)^{1/2}.$
Introducing $v=\left( \begin{array}{c} \langle D_x\rangle u \\ u^\prime
\end{array}\right),$ the matrix form of this equation is given by
$$ v^\prime=\left( \begin{array}{cc} 0 & \langle D_x\rangle \\
-\langle D_x\rangle^{-1}c  & -b \end{array}\right) v.$$
Let $b_1$ and $c_2$ be principal symbols of $b$ and $c$. The equation is
hyperbolic if $b_1^2\geq 4c_2$, with multiple roots at $b_1^2=4c_2.$
Assume that $b_1^2-4c_2=\mu^2$, with $\mu\in S^1$ being a symbol of order
one. Then characteristics $a_1, a_2$ satisfy
$\{a_1,a_2\}=\frac12\{b_1,\mu\},$ which may vanish.

Let the Cauchy data be $u(0)=f_0, u^\prime(0)=f_1.$ Let
$\alpha=(n-1)|1/p-1/2|$ and
$1<p<\infty.$ If $\mu$ is elliptic, it is known that if $f_j\in
L^p_{\alpha-j}$, then
$u(t,\cdot)\in L^p.$
Our result of Theorem \ref{th:lp}
will imply that the same is true for any $\mu\in S^1.$

Note that in this example we may requre only microlocal diagonalization
with the same conclusion.

Let us now discuss the propagation of singularities for operators $Q_l$.
This result is essentially a reformulation of Rozenblum's result in the
case of finite geometry under Condition C.
It is clear (also from multi-phase analysis) that singularities propagate
along broken Hamiltonian flows. Let
$$J=\{j_1,\ldots,j_{l+1}\}, \; 1\leq j_k\leq m,\; j_k\not=j_{k+1}.$$
Let $\Phi_J(t,x,\xi)$
be the corresponding broken Hamiltonian flow. It means that points
follow bicharacteristics of $a_{j_1}$ until meeting the characteristic
of $a_{j_2}$, and then continue along the bicharacteristic of $a_{j_2}$,
etc.
Note that singularities may accumulate if wave front sets for different
broken trajectories project to the same point of $X$.

We can write
$$Q_l=\int_\Delta I(\bar{t}) d\bar{t},$$
where $\bar{t}=(t_1,\ldots,t_l)$,
$\Delta=\{0\leq t_l\leq t_{l-1}\leq \ldots\leq t_1\leq t\}$
is a symplex in $\Rone^l$ and
$I(\bar{t})=Z(t_1)\circ\ldots\circ Z(t_l).$
It is possible to treat it as
a standard Fourier integral operator
with the change of variables $\bar{t}=\zeta|\xi|^{-1}.$
Let $K$ be a cone in $\Rone^N=\Rone^{n+l}.$
Let
$$Iu(x)=\int_K\int_Y e^{i\varphi(x,y,\theta)} a(x,y,\theta) u(y) dy d\theta$$
be a Fourier integral operator with integration over the cone $K$ with respect
to $\theta.$
Let $K_j$ be $K$ or a face of $K$. Let $\varphi_j(x,y,\theta_j)=\varphi|_{K_j},
\theta_j\in K_j.$ Let $\Lambda_j\subset T^*X\times T^*X$ be a
Lagrangian manifold with boundary:
$$\Lambda_j=\{(x,\frac{\partial\varphi_j}{\partial x},y,
-\frac{\partial\varphi_j}{\partial y}): \frac{\partial\varphi_j}{\partial\theta_j}=0\}.$$
For $G\subset T^*Y$,
let $\Lambda_j(G)=\{z\in T^*X: \exists\zeta\in G: (z,\zeta)\in\Lambda_j\}.$
Then we have the following statement on the propagation of singularities.

\begin{theorem}
Let $u\in {\mathcal{D}}^\prime(Y).$ Then
$WF(Iu)\subset \cup_j \Lambda_j(WF(u)).$
\label{th:singularities}
\end{theorem}
The proof is standard and follows H\"ormander \cite{H}.

From this, we can deduce first and second terms of the spectral
asymptotic of operator $P$.
Let us call $T$ a period of symbol $A(x,\xi)$ if there exists
$J$ such that $j_1=j_{l+1}$, and the trajectory of $\Phi_J$ is closed:
$\Phi_J(T,x,\xi)=(x,\xi).$
Then we have the following extension of well-known results of
H\"ormander \cite{Hormander-sp}, Duistermaat--Guillemin \cite{DG},
Safarov--Vassiliev \cite{SV}, and Rozenblum \cite{Rozenblum}.

\begin{theorem}
Assume that $X$ is compact and assume that Condition C is
satisfied. Let $D$ be the
set of $(x,\xi)\in T^*X$ such that there exist $T$ and $J$ such that
$\Phi_J(T,x,\xi)=(x,\xi).$
Assume that the measure of $D$ is zero.
Then for the spectrum of $P$ the following
Weyl formula holds:
$$N(\lambda)=\sharp\{j: \lambda_j<\lambda\}=c_n \lambda^n+c^\prime_n \lambda^{n-1}
+o(\lambda^{n-1}),$$
where $\lambda_j$ are eigenvalues of $P$.
\label{th:spectral}
\end{theorem}
Proof of the $L^p$ estimates will be based on Theorem \ref{th:regularity},
which we regard as an independent result on regularity of solutions of
partial differential equations. Theorem \ref{th:eps} concerns the measure
of the set where a function is small given some information on the
multiplicity of its roots. It will play a crucial role in the proof of
the smoothing property in Theorem \ref{th:q}.

\section{Regularity of solutions}

In this section we will present a principle governing solutions of
first order systems.
Let $Z(t)\in\Lcal(C_0^\infty(X),\Dcal^\prime(X))$, $t\in [0,T]$,
be a time dependent family of operators.
Let $W_0, W_1$ and $W$ be linear subspaces of $\Dcal^\prime(X)$
such that $W_0, W_1\hookrightarrow W.$
We will make different choices of these spaces in the future, dependent
on the structure of operators $Z(t).$
Let us consider the Cauchy problem for $u=u(t,x)$:
\begin{equation}
\left\{
\begin{array}{ll}
  u^\prime-Z(t)u = r, & r(t)\in W_0, \\
  u(0)\in W_1. &  \\
\end{array}
\right.
\label{eq:eqpr1}
\end{equation}
One is often interested in the following question. If the right hand side
and Cauchy data satisfy $r(t)\in W_0$ and $u(0,\cdot)\in W_1$,
when do fixed time solutions $u(t,\cdot)$ of (\ref{eq:eqpr1})
belong to $W$?
In general, some loss of regularity is possible in problem
(\ref{eq:eqpr1}) even if $Z(t)$ are very good.
So we will think of $W_0$ being the smallest,
$W_1$ an intermediate, and $W$ the largest among these spaces.
The following theorem says that if operators $Z(t)$ have some
structure, and solutions of Cauchy problem (\ref{eq:eqpr1}) with
zero Cauchy data are in $W$, so will be solutions with Cauchy
data from some sufficiently large space $W_1.$

\begin{theorem}
Let $W_0, W_1\hookrightarrow W$ be linear subspaces of
$\Dcal^\prime(X)$. \\
Let $Z(t)\in \Lcal(C_0^\infty(X),\Dcal^\prime(X))$, $t\in [0,T]$.
Assume that
\begin{itemize}
\item[{\rm (i)}]{\rm (Boundedness)} $Z$ extends to an operator
in $L^\infty([0,T],\Lcal(L^2(X),L^2(X)))$, and
$Z(t)$ extend to continuous linear
operators from $W_1$ to $W$, for all $t\in [0,T].$
\item[{\rm (ii)}]{\rm (Calculus)}
$Z(t_1)\circ\cdots\circ Z(t_l):W_1\to W$ are continuous
for all $l$ and for all $t_1,\ldots,t_l\in [0,T]$.
\item[{\rm (iii)}]{\rm (Smoothing)} There exists $l$ such that
$$Z(t)\int_0^t\int_0^{t_1}\cdots\int_0^{t_{l-1}} Z(t_1)\circ
\cdots\circ Z(t_l) dt_l\cdots dt_1$$
is continuous from $W_1$ to $W_0,$ for all $t\in [0,T]$.
\item[{\rm (iv)}]{\rm (Zero Cauchy data)}
Solutions $v=v(t,x)$ of the Cauchy problem
  \begin{equation}
     \left\{ \begin{array}{ll}
         v^\prime-Z(t)v=r, & r(t)\in W_0,
          \\
         v(0)=0,
      \end{array} \right.
      \label{eq:cc}
    \end{equation}
  satisfy $v(t,\cdot)\in W$ for $t\in [0,T].$
\end{itemize}
Then the solution $u=u(t,x)$ of the Cauchy problem
 \begin{equation}
     \left\{ \begin{array}{ll}
         u^\prime-Z(t)u=r, & r(t)\in W_0,
          \\
         u(0)\in W_1,
      \end{array} \right.
      \label{eq:cauchy1}
    \end{equation}
  satisfies $u(t,\cdot)\in W$ for all $t\in [0,T]$.

  Moreover, if $W_0, W_1, W$ are normed spaces and if
  solutions $v(t,\cdot)$ to \eqref{eq:cc} in {\rm (iv)}
  satisfy $\| v(t,\cdot)\|_W\leq C\| r(t)\|_{W_0}$ for all
  $t\in [0,T]$, then also
  $$\| u(t,\cdot)\|_W\leq C(\|u(0)\|_{W_1}+\|r(t)\|_{W_0}),$$
  for all $t\in [0,T]$.
  \label{th:regularity}
\end{theorem}
Conditions (i) and (ii) ensure that operators $Z(t)$ have
some structure. Indeed, if $W_1\subset W$ is different from $W$,
(ii) does not follow from (i). In our typical applications,
$Z(t)$ will be time dependent pseudo-differential or Fourier
integral operators, and compositions in (ii) are essentially
of the form of a single operator $Z(t)$. Condition (iii) is
natural from the point of view of harmonic analysis, since
integration with respect to a parameter often brings additional
regularity. Condition (iv) ensures that solutions with zero
Cauchy data and regular right hand side are also sufficiently regular.
\begin{proof}
Let $U(t)$ be an operator solving the Cauchy problem
\begin{equation}\label{eq:cauchyop}
\left\{
\begin{array}{ll}
  U^\prime-Z(t)U = R(t), & R(t)\in \Lcal(W_1,W_0), \\
  U(0)=I. &  \\
\end{array}
\right.
\end{equation}
Let $U_0(t)$ be some partial solution to the problem
\[
\left\{
\begin{array}{l}
 U_0^\prime-Z(t)U_0(t)=R(t), \\
 U_0(0)=0.
\end{array}
\right.
\]
Then the solution $U$ of problem (\ref{eq:cauchyop}) satisfies
\begin{equation}
U(t)=U_0(t)+I+\int_0^t Z(t_1) dt_1 + \int_0^t\int_0^{t_1} Z(t_1)
Z(t_2) dt_2 dt_1 + \ldots
\label{eq:eqseries}
\end{equation}
The convergence of this series can be understood in $L^2$. Indeed,
because of assumption (i), the term of this series with $k$
integrals can be estimated by $t^k \sup_t ||Z(t)||^k_{L^2\to
L^2}/k!$ From this it also follows that $U(t)$ is a solution of
(\ref{eq:cauchyop}) in $L^2$.
Let us now define
\begin{align*}
S_N(t)=I+\int_0^t Z(t_1) dt_1 + \int_0^t\int_0^{t_1} Z(t_1) Z(t_2)
dt_2 dt_1 + \ldots
\\
+ \int_0^t\int_0^{t_1}\ldots \int_0^{t_{N-1}} Z(t_1) Z(t_2)\ldots
Z(t_N) dt_N\ldots dt_2 dt_1.
\end{align*}
Let $V(t)=U(t)-S_N(t)$, it is equal to $U_0(t)$ plus the remainder
of the series (\ref{eq:eqseries}). Then we have
\begin{eqnarray}
V^\prime(t)-Z(t)V(t) = \nonumber
\\
(U-S_N)^\prime(t)-Z(t)(U-S_N)(t) = \nonumber
\\
(U^\prime-ZU)(t)-(S_N^\prime-ZS_N)(t) = \nonumber
\\
R(t)-Z(t)\int_0^t\int_0^{t_2}\ldots \int_0^{t_{N-1}} Z(t_2)\ldots
Z(t_N) dt_2\ldots dt_N.
\label{eq:aux1}
\end{eqnarray}
Choosing $N=l$, from assumption (iii) of the theorem the second term
is continuous from $W_1$ to $W_0$. Since also
$R(t)\in\Lcal(W_1,W_0)$, it follows that the right hand side is a
continuous linear operator from $W_1$ to $W_0$.

Let $w=u(0)\in W_1$ be the Cauchy data for \eqref{eq:cauchy1}.
If we denote by $\rho(t)$
the value of the operator in the last line of \eqref{eq:aux1}
at $w$, we will have $\rho(t)\in W_0.$
The value of $V(0)$ is
\[
V(0)=U(0)-S_N(0)=0.
\]
It follows now that $V(t)w$ solves Cauchy problem (\ref{eq:cc}),
so it belongs to $W$ by assumption (iv).
Since $S_N(t)$ is continuous from $W_1$ to $W$ by assumption (ii),
and $V(t)w=U(t)w-S_N(t)w$ is in $W$, be obtain
$u(t,\cdot)=U(t)w\in W$.

Moreover, suppose that we also have the estimate
$\|v(t,\cdot)\|_W\leq C\|\rho(t)\|_{W_0}$ in (iv).
Then we also have
$$
\| u(t,\cdot)\|_W\leq \|V(t)w\|_W+\|S_N(t)w\|_W\leq
C\|\rho(t)\|_{W_0}+C\|w\|_{W_1}\leq C(\|w\|_{W_1}+\|R(t)w\|_{W_0}).
$$
\end{proof}

Later we will need this in the case of $Z(t)$ being Fourier
integral operators of order zero. However, let us point out several
applications to other cases of pseudo-differential and
Fourier integral operators.
In these cases we will make different appropriate choices of
spaces $W_0, W_1, W$. Moreover, in Corollaries
\ref{cor:c0} and \ref{cor:c1} we will assume that the corresponding
non-homogeneous Cauchy problems with zero Cauchy data have unique
solutions. This is a natural assumption if $Z(t)$ behave sufficiently
well with respect to $t$ since we are working in
subspaces of $L^2$.

\begin{corollary}\label{cor:c0}
Let $1<p<\infty$. Let $Z(t)\in \Psi^0$, $t\in [0,T]$, be a family of
pseudo-differential operators of order zero with amplitudes
compactly supported in $x,y$, uniformly in $t$. Suppose
that $Z\in L^\infty([0,T],\Lcal(L^2,L^2))\cap
L^\infty([0,T],\Lcal(L^p,L^p))$. Then the
solution $u=u(t,x)$ of the Cauchy problem
\begin{equation}\label{eq:cauchyz}
\left\{
\begin{array}{ll}
  u^\prime-Z(t)u = r, & r(t)\in L^p,\ t>0, \\
  u(0)\in L^p, &  \\
\end{array}
\right.
\end{equation}
satisfies $u(t,\cdot), u^\prime(t,\cdot)\in L^p$, for all
$t\in (0,T].$ Moreover, we have an estimate
$$\sup_{t\in [0,T]} \|u(t,\cdot)\|_{L^p}\leq C\|u(0)\|_{L^p}.$$
\end{corollary}
If amplitudes of $Z(t)$ are not compactly supported with respect to
$x, y$, we have a similar local statement for compactly supported
Cauchy data.
Note also that pseudo-differential operators of order zero are locally
bounded in $L^2$ and $L^p$, for all $1<p<\infty$. Conditions
$Z\in L^\infty([0,T],\Lcal(L^2,L^2))$ and
$Z\in L^\infty([0,T],\Lcal(L^p,L^p))$ simply mean that we have some
control on their norms, i.e. there exist a constant $C$ such that
$$
\|Z(t)\|_{L^2\to L^2}\leq C,\
\|Z(t)\|_{L^p\to L^p}\leq C,\ \forall t\in [0,T].
$$
\begin{proof}
Let us choose $W=W_1=W_0=L^p_{comp}$.
Let us check conditions of Theorem \ref{th:regularity}.
Properties (i) and (ii) follow from regularity properties of
pseudo-differential operators of order zero and our assumptions.
Property (iii) also
holds because $Z(t)$ are locally bounded in $L^p$.
Property (iv) is a consequence of Duhamel's principle and is similar
to the one in Corollary \ref{cor:c1}. Picard series is convergent
in $L^p$ provided that operators norms $\|Z(t)\|_{L^p\to L^p}$
are uniformly bounded for $t\in [0,T].$ Norm estimate follows
from this as well.

We can see that $u^\prime(t,\cdot)\in W$ from $u^\prime=Z(t)u+r$ and
from the continuity of $Z(t)$ in $W$.
\end{proof}

We will now apply Theorem \ref{th:regularity} in the case of
$Z(t)$ being Fourier integral operators. While our case
\eqref{eq:eq4} corresponds
to $Z(t)$ being operators of order zero, the crucial smoothing
property (iii) will follow from the fact that operators $Z(t)$ have
a special structure. For general Fourier integral operators
$Z(t)$ without structure, we need to assume that they are of negative orders.
This
is for example the case when the zero order term $B$ in
Theorem \ref{th:lp} is actually a pseudo-differential
operator of some negative order.

\begin{corollary}\label{cor:c1}
Let $1<p<\infty$, $\epsilon>0$, and $\alpha=(n-1)|1/p-1/2|.$
Let $Z(t)\in I^{-\epsilon}$, $t\in [0,T]$, be a family of non-degenerate
Fourier integral operator of order $-\epsilon$ with amplitudes compactly
supported in $x,y$, uniformly in $t$. Suppose that operators $Z(t)$ can
be composed and that
$Z\in L^\infty([0,T],\Lcal(H^s,H^s))$, for some $s>(pn-2n)/2p$ when
$p>2$ and $s=0$ when $p\leq 2$.
Then the solution $u=u(t,x)$ of
the Cauchy problem
\begin{equation}\label{eq:cauchyzf}
\left\{
\begin{array}{ll}
  u^\prime-Z(t)u = r, & r(t)\in H^s,\ t>0, \\
  u(0)\in (L^p_\alpha)\cap L^2, &  \\
\end{array}
\right.
\end{equation}
satisfies $u(t,\cdot)\in L^p$, for all $t\in
(0,T].$ Moreover,
$$ \sup_{t\in [0,T]}\|u(t,\cdot)\|_{L^p}\leq C\|u(0)\|_{L^p_\alpha}.$$
\end{corollary}
Note that operators $Z(t)$ are locally bounded in $H^s$, so
assumption \\
$Z\in L^\infty([0,T],\Lcal(H^s,H^s))$ simply means
that $\|Z(t)\|_{H^s\to H^s}\leq C$ for all $t\in [0,T]$.
\begin{proof}
Let $W=L^2\cap L^p_{comp}$, $W_1=(L^p_\alpha)_{comp}$, and
$W_0=H^s_{comp}\subset W$.
Let us check conditions of Theorem \ref{th:regularity}. Condition
(i) follows from the fact that non-degenerate Fourier integral
operators of order 0 are bounded from $(L^p_\alpha)_{comp}$ to
$L^p_{loc}$.
Condition (ii) follows from the calculus of
non-degenerate Fourier integral operators, since we assumed that
compositions of $Z(t)$ are again non-degenerate Fourier integral
operators.
Smoothing condition (iii)
for large $l$ follows again from the calculus, since operators
$Z(t)$ are of order $-\epsilon$.

Finally, let us show that solutions of $v^\prime-Z(t)v=r(t)$,
$r(t)\in H^s_{comp}$, with zero Cauchy data $v(0)=0$,
satisfy $v(t,\cdot)\in L^2\cap L^p$.
In fact, we will show that $v(t,\cdot)\in H^s\subset
L^p\cap L^2$.

From the uniqueness of the solution of this problem
it follows that we can use Duhamel's principle to write
\begin{equation}\label{eq:repr}
v(t,x)=\int_0^t E(t,s) r(s,x) ds,
\end{equation}
where $E(t,s)$ is the propagator of
\[
\left\{
\begin{array}{l}
  (\partial_t-Z(t))E(t,s) = 0, \\
  E(t,s)|_{t=s}=I.  \\
\end{array}
\right.
\]
Picard series for this problem gives the asymptotic expansion of
$E(t,s)$, in particular implying that $E(t,s)$ is bounded in $L^2$ and in
$H^s$ provided operator norms $\|Z(t)\|_{H^s\to H^s}$ are uniformly
bounded for $t\in [0,T]$. From (\ref{eq:repr}) it follows
that $v(t,\cdot)\in H^s$.  Moreover, since $r(s,\cdot)\in H^s$, we
also get an estimate
$$
\|u(t,\cdot)\|_{H^s}\leq C\sup_{\tau\in [0,T]}\|r(\tau,\cdot)\|_{H^s},
$$
implying the estimate in Corollary \ref{cor:c1}.
\end{proof}

\begin{proof}[Proof of Theorem \ref{th:lp}]
As we have already mentioned, by Lemma \ref{l:diag}
we can assume that characteristics of
$A$ are correctly defined on $T^*X$.
Since $A$ is diagonalizable, we can write
\begin{eqnarray}\label{eq:diag}
P(t,x,D_x) & = & \bigoplus a_j(x,D_x)+B(t,x,D_x), \\
B(t,x,D_x) & = & (B_{jk}(t,x,D_x))_{1\leq
j,k\leq m}, \ B_{jk}\in C^\infty([0,T], \Psi^0). \nonumber
\end{eqnarray}
Some of $a_j$'s may be identically equal to each other.
We can renumber $a_j$'s into $r$ groups (possibly of size one) of
equal characteristics. These are the eigenvalues of the matrix
$A(x,\xi)$ counted with multiplicity. Thus, we have
$1=k_1<k_2<\ldots<k_r=n+1,$
and $a_{k_i}\equiv\ldots\equiv a_{k_{i+1}-1}\not\equiv a_k$, for
$k<k_i$ or
$k\geq k_{i+1}$. This means that we have a group of the same roots
$a_1,\ldots,a_{k_2-1}$, etc., while roots from different groups
are not identically the same.
Therefore, this is a decomposition of the first order principal part
into a block-diagonal form with the same roots in each block, with possible
equality of roots in different blocks at some points.
So we can write
\begin{equation}\label{eq:diag1}
P(t,x,D_x)=\diag(\tilde{a}_1,\ldots,\tilde{a}_r)+B(t,x,D_x),
\end{equation}
where $\tilde{a}_i=\diag (a_{k_i},\ldots,a_{k_{i+1}-1})$ are diagonal
matrices with equal roots at the diagonal. Let us set
\[
\tilde{A}_i=\tilde{a}_i+(B_{\mu\nu})_{k_i\leq\mu,\nu\leq k_{i+1}-1},
\]
so that
$$
P=\tilde{A}+B=\diag (\tilde{A}_1,\ldots,\tilde{A}_r)+B.
$$
Note that in the last equality we can
assume $B_{\mu\nu}=0$ for $k_i\leq\mu,\nu\leq k_{i+1}-1$ if we add these
components to the corresponding components of $\tilde{A}$.
Let $U(t)=\exp(-i\tilde{A}t)V(t)$. This is well defined in view of,
for example, \cite[VIII]{Taylor}. Then
\[
V^\prime=Z(t) V\equiv -i e^{i\tilde{A}t} B e^{-i\tilde{A}t} V, \;\; V(0)=I.
\]
Now we will apply Theorem \ref{th:regularity} with
$Z(t)=-i e^{i\tilde{A}t} B e^{-i\tilde{A}t}$. Let us choose
$W=L^2_{comp}\cap L^p$, $W_1=(L^p_\alpha)_{comp}$, and
$W_0=H^s_{comp}$ with $s>(pn-2n)/2p$ for $p>2$ and
$s=0$ for $1<p\leq 2$. Conditions (i) and (ii) follow from the
calculus and regularity properties of non-degenerate
Fourier integral operators
of order zero. Smoothing condition (iii) follows from Theorem
\ref{th:q}. For condition (iv) we can use Duhamel's principle
similar to the proof of Corollary \ref{cor:c1}. Thus, Theorem
\ref{th:regularity} implies that propagator $V(t)$ is continuous
from $L^p_\alpha$ to $L^p$. Operator $U(t)$ is a composition
of $V(t)$ with a non-degenerate Fourier integral operator
$\exp(-i\tilde{A}t)$, so $U(t)$ is given as a sum of a smoothing
series obtained by the multiplication of Picard series for
$V(t)$ with $\exp(-i\tilde{A}t)$. Using the calculus of Fourier
integral operators in each term of the series and its smoothing
property we can repeat the argument of Theorem \ref{th:regularity}
in this case to see that $u(t,\cdot)\in L^p$ with an estimate
for its norm.
\end{proof}
Note that if $B$ in Theorem \ref{th:lp} is a pseudo-differential
operator of negative order, $B\in\Psi^\mu$, for some $\mu<0$,
the proof is simpler because we do not have to use Theorem
\ref{th:q} to prove condition (iii) of Theorem \ref{th:regularity}.
Instead, we can use directly Corollary \ref{cor:c1} to obtain the
smoothing property (iii).

\begin{proof}[Proof of Theorem \ref{th:inv}]
Let $L$ be as in \eqref{eq:inv1} and let
$u$ be the solution of
\be
Lu=f,\ D_t^j u(0,x)=g_j(x),\ 0\leq j\leq m-1.
\label{eq:meq}
\ee
Let
$$
U=(u,\partial_1 u,\partial_2 u,\ldots,
\partial_1\partial_2 u,\partial_2\partial_1 u,\ldots,
\partial^J u,\ldots)^T,
$$
where $J=\{j_1,\ldots,j_k\}$ is a permutation of some elements
of $\{1,\ldots,m\}$, $|J|=k\leq m-1$. 
Vector $U$ has $1+\sum_{j=1}^{m-1}m!/j!$
components. Here we can write
$\partial^J=D_t^k+\sum_{j=0}^{k-1} c_j^J(t,x,D_x)D_t^j$,
where $c_j^J(t)\in \Psi^{k-j}$. We set $|J|=k$. It was
shown by Morimoto in \cite{Morimoto} that $U$ solves
the system
\be
D_t U + AU+BU=F,\ U(0,x)=G(x),
\label{eq:msys}
\ee
where $F=(0,\ldots,0,f,\ldots,f)$ and
$G=(g_0,\ldots,g_{|J|}+\sum_{j=0}^{|J|-1} c_j^J g_j,\ldots),$
$A$ is a diagonal matrix with $\lambda_j$'s at the diagonal
and $B$ is a matrix of pseudo-differential operators of
order zero. Matrix $B$ has some operators in the last row,
zeros, and $-1$ above the diagonal. 
If $\lambda_j$ satisfy Condition C, Theorem
\ref{th:lp} implies
$\|U(t,\cdot)\|_{L^p}\leq C\|G\|_{L^p_\alpha}.$
Since $c_j^J(t)\in\Psi^{k-j}$, we get
$\|G\|_{L^p_\alpha}\leq
C \sum_{j=0}^{m-1}\|g_j\|_{L^p_{\alpha+m-1-j}}$,
which implies the estimate of the Theorem.
\end{proof}

\section{Estimates for Picard series}

In this Section we will prove Theorem \ref{th:q} on the smoothing
properties of terms $Q_l$ of the Picard series \eqref{eq:series}.

Let $A_j\in\Psi^1$, $j=1,\ldots,r$, be elliptic pseudo-differential
operators of order one. Let $a_j(x,\xi)$ denote their
principle symbols. We can assume that there are no idential symbols
among these $a_j$'s.
Let
\[H(\bar{t})=e^{iA_{j_1}t_1}e^{iA_{j_2}(t_2-t_1)}\cdots
e^{-iA_{j_{l+1}}t_l},\]
where $1 \leq j_k \leq r$, $j_k\neq j_{k+1}$, $k=1,\ldots,l+1$, and
$\bar{t}=(t_1...t_l)$.
Let us define
\[Q=\int_0^t\int_0^{t_1}...\int_0^{t_{l-1}}B(\bar{t})H(\bar{t})dt_l...dt_1,
\]
where $B(\bar{t})\in\Psi^0$ is a pseudo-differential operator of order zero
smoothly dependent on $\bar{t}$.
Such operators $Q$ appear in the Picard series \eqref{eq:series},
\eqref{eq:ql}.
In this section we will give a detailed description of operator $Q$
in order to prove that it is a smoothing operator when $l$ is
sufficiently large. First of all let us note that $H$ is a
Fourier integral operator and due to the
theorem on compositions of Fourier integral operators
the canonical relation $\Lambda^{\bar{t}} \subset T^*X\times T^*X$ of
$H(\bar{t})$ is given by
\[\Lambda^{\bar{t}}=\{(x,p,y,\xi): (x,p)=\Psi^{\bar{t}}(y,\xi)\},
\]
where
$\Psi^{\bar{t}}=\Phi_{j_1}^{t_1}\circ\cdots
\Phi_{j_l}^{t_l-t_{l-1}}\circ
\Phi_{j_l+1}^{-t_l}$ and
$\Phi_j^t$ is the Hamiltonian flow defined by $a_j$.

It can be easily checked that  $H$ is a solution operator for
the system of equations
\be
\frac{\partial H}{\partial
t_k}=iT_k(t_1,...,t_k)H, \ k=1,\ldots,l,
\label{hyp}
\ee
where $T_k\in\Psi^1$ is a pseudo-differential operator of order one.
In view of Egorov's theorem its principle symbol is equal to
\be
T_k^0(t_1,...,t_k,x,\xi)=(a_{j_k}-a_{j_{k+1}})\circ
\Phi_{j_1}^{t_1}\circ\cdots\circ\Phi_{j_k}^{t_k-t_{k+1}}(x,\xi),\
\ (x,\xi)\in T^*X
\label{prin}
\ee
for all $k=1,\ldots,l$.
Let us construct a phase function
 $\varphi(\bar{t},x,y,\xi)$ which defines the
operator $H(\bar{t})$ for small $\bar{t}$. We will look for it in
the form $\varphi(\bar{t},x,y,\xi)=\psi(\bar{t},x,\xi)-y\cdot\xi$.
It follows from \eqref{hyp} and \eqref{prin} that $\psi$ satisfies
a system of Hamilton-Jacobi equations
\be
\frac{\partial\psi}{\partial t_k}=T_k^0(\bar{t},x,
\frac{\partial\psi}{\partial x}),\
\psi(0,x,\xi)=x\cdot\xi.
\label{hj}
\ee
In \cite{Rozenblum} it was
checked that Frobeneus conditions for system
\eqref{hj} are satisifed.
Solving this system we obtain a non-degenerate phase
function. This phase function defines a Lagrangian manifold
$\Lambda^{\bar{t}}$, so that we have
\be
(x,\frac{\partial\psi}{\partial
x})=\Psi^{\bar{t}}(y,\xi)=(x^{\bar{t}}(y,\xi),p^{\bar{t}}(y,\xi)),
\ y=\frac{\partial\psi}{\partial \xi}.
\label{lagr}
\ee
Now we are going to investigate
the smoothing properties of operator $Q$.
We can write $Q$ as
$$
Qu(x)=\int_\Delta\int_\Rn\int_\Rn e^{i\varphi(\bar{t},x,y,\xi)}
b(\bar{t},x,y,\xi) u(y) dy d\xi d\bar{t},
$$
where $\varphi$ satisfies \eqref{lagr} and $b$ is an
amplitude of order zero, which we may assume to be compactly
supported with respect to $x$ and $y$. Here $\Delta$ is a
symplex $\{0\leq t_l\leq t_{l-1}\leq\ldots\leq t_1\leq t\}.$

Let $\chi\in C_0^\infty$
be a cut-off function such that $\chi(\tau)=1$ for $|\tau|<1$ and
$\chi(\tau)=0$ for $|\tau|>2$. Operator $Q$ can be decomposed as
$Q=R_1+R_2+R_3$, where
$$
R_ju(x)=\int_\Delta\int_\Rn\int_\Rn e^{i\varphi(\bar{t},x,y,\xi)}
\mu_j(\epsilon,\bar{t},x,y,\xi)u(y)dyd\xi
d\bar{t}, \ j=1,2,3,
$$
where
\begin{eqnarray*}
\mu_1(\epsilon,\bar{t},x,y,\xi) & = &
\left(1-\chi(\epsilon^{-1}|\frac{\partial
\varphi}{\partial \xi}|)\right)b(\bar{t},x,y,\xi), \\
\mu_2(\epsilon,\bar{t},x,y,\xi)& = &
\chi(\epsilon^{-1}|\frac{\partial
\varphi}{\partial \xi}|)\left(1-\chi(\epsilon^{-1}|\frac{\partial
\varphi}{\partial \bar{t}}||\xi|^{-1})\right)b(\bar{t},x,y,\xi), \\
\mu_3(\epsilon,\bar{t},x,y,\xi) & = &
\chi(\epsilon^{-1}|\frac{\partial
\varphi}{\partial \xi}|)\chi(\epsilon^{-1}|\frac{\partial
\varphi}{\partial \bar{t}}||\xi|^{-1})b(\bar{t},x,y,\xi).
\end{eqnarray*}
Let us first consider operator $R_1$. On the support of $\mu_1$
we have the estimate
$$
\left|\frac{\partial \varphi}{\partial \xi}\right|\geq \epsilon.
$$
Therefore, there exists an operator
$L(\frac{\partial}{\partial \xi})$ of order 1, such that
$Le^{i\varphi}=e^{i\varphi}$, with coefficients estimated by
$C\epsilon^{-1}$ on the support of the amplitude of $\mu_1$.
When integrating by parts with $L$ there may appear an additional
factor of $\epsilon^{-1}$ when differentiating $\chi$.
Integrating by parts $p$ times with this operator $L$ we obtain
an operator with an amplitude of order
$-p$ and coefficients that can be estimated by $\epsilon^{-2p}$.
From Lemma \ref{rozen:modify} with $q=-p$ we obtain the following
estimate
\be
\|R_1\|_{p-n-1}\leq
C\epsilon^{-3p+n+1}.
\label{r1}
\ee
The same procedure with
integrating by parts with respect to $\xi$ can not be applied to
$R_2$. But here there is a possibility to integrate by parts with
respect to $\bar{t}$. Indeed, there exists an operator
$M(\frac{\partial}{\partial \bar{t}})$, such that
$Me^{i\varphi}=e^{i\varphi}$, with coefficients not greater than
$\epsilon^{-1}|\xi|^{-1}$ on the support of $\mu_2$. Integrating by
parts with $M$ we obtain an operator with an amplitude of order $-1$ with
coefficients that can be estimated by $\epsilon^{-2}$, where
another $\epsilon^{-1}$ may appear from differentiating $\chi$.
The boundary
integrals have the same form as $Q$ but they have amplitudes of order
$-1$ and depend on not less than $l-2$ time variables.
The reason for possibly losing two variables is that after restriction
to the boundary, say $t_2=t_1$, it may happen that $a_{j_1}$ and $a_{j_3}$
are the same roots.
It follows that
we can apply such procedure $[l/2]$ times. As a result we
obtain operators of order $-[l/2]$ with coefficients not greater
than $C\epsilon^{-2[l/2]}$. Then in view of Lemma \ref{rozen:modify}
for $[l/2]>n+1$ we have
\be
\|R_2\|_{[l/2]-n-1}\leq C\epsilon^{-3[l/2]+n+1}.
\label{r2}
\ee
Let us now consider
the last integral $R_3$. It is not possible to apply procedures with
integrating by parts as above either with respect to $\xi$ or with
respect to $\bar{t}$. But in this case it turns out that the
support of amplitude $\mu_3$ is small. The singular support of the
integral kernel of $Q$ may be very irregular in this case, so a more
delicate analysis is necessary to show the smoothing properties of $R_3$.
First we will show that
 \be
 |T_j^0(\bar{t},x^{\bar{t}}(y,\xi),p^{\bar{t}}(y,\xi))|\leq C\epsilon|\xi|,
 j=1,\ldots,l,
 \label{res}
 \ee
on the support of $\mu_3$. We notice that $\mu_3$ differs from $0$ only if
\be
|\frac{\partial\varphi}{\partial\xi}|\leq 2\epsilon , \ \ \
|\frac{\partial\varphi}{\partial\bar{t}}|\leq 2\epsilon|\xi|.
 \label{obl}
\ee
 It follows from \eqref{hj} and \eqref{obl} that
\be
|T_j^0(\bar{t},x,\frac{\partial \psi}{\partial x})|\leq C\epsilon |\xi|,
j=1,\ldots,l.
\label{1}
\ee
  Because of homogeneity of $T_j^0$ with respect to $\xi$ we also
  have the following trivial estimates
$$
|\partial_x T_j^0(\bar{t},x,\frac{\partial \psi}{\partial
 x})|\leq C|\xi|, \ \ \
 |\partial_\xi T_j^0(\bar{t},x,\frac{\partial \psi}{\partial
 x})|\leq C,\ \ j=1,\ldots,l.
$$
 Consequently, we obtain
\begin{eqnarray*}
|T_j^0(\bar{t},x,\frac{\partial \psi}{\partial
 x})-T_j^0(\bar{t},x^{\bar{t}}(y,\xi),p^{\bar{t}}(y,\xi))| & \\
 \leq
 C(|\xi||x-x^{\bar{t}}(y,\xi)|+|\frac{\partial \psi}{\partial
 x}-p^{\bar{t}}(y,\xi)|),\ \  j=1,\ldots,l. & \nonumber
 \label{2}
\end{eqnarray*}
Since $\frac{\partial\varphi}{\partial\xi}(\bar{t},x^{\bar{t}}(y,\xi),y,\xi)=0$,
we get
$$
 \frac{\partial\varphi}{\partial\xi}(\bar{t},x,y,\xi)=
 \frac{\partial\varphi}{\partial\xi}(\bar{t},x,y,\xi)-
 \frac{\partial\varphi}{\partial\xi}(\bar{t},x^{\bar{t}}(y,\xi),y,\xi)=
 \frac{\partial^2 \varphi}{\partial x\partial\xi}(\bar{t},x^*,y,\xi)
 (x-x^{\bar{t}}(y,\xi)),
$$
for some points $x^*$.
Since $|\frac{\partial^2\varphi}{\partial x \partial \xi}|\neq 0$ for small
 $\bar{t}$, we obtain
\be
|x-x^{\bar{t}}(y,\xi)|\leq C \left|\frac{\partial\varphi}{\partial
 \xi}(\bar{t},x,y,\xi)\right|.
 \label{x}
\ee
From the properties of Lagrangian manifold $\Lambda^{\bar{t}}$ in \eqref{lagr}
we also obtain
\be
|\frac{\partial \psi}{\partial
 x}(\bar{t},x,\xi)-p^{\bar{t}}(y,\xi)|=|p^{\bar{t}}(\frac{\partial \psi}{\partial
 \xi},\xi)-p^{\bar{t}}(y,\xi)|\leq C
 |\frac{\partial \varphi}{\partial \xi}(\bar{t},x,y,\xi)||\xi|,
 \ \ \forall (x,y,\xi).
   \label{p}
\ee
Finally, taking into consideration \eqref{obl}-\eqref{p} we obtain \eqref{res}
on the support of $\mu_3$.

Estimate \eqref{res} shows us that the support of the amplitude $\mu_3$ of
operator $R_3$ is contained in
$$
\Xi_1=\{(\bar{t},x,y,\xi): x\in X,\ \xi\not=0,\
(\bar{t},y,\frac{\xi}{|\xi|})\in \Xi \},
$$
where
$$
\Xi=\{(\bar{t},y,\frac{\xi}{|\xi|}):\
|T_j^0(\bar{t},x^{\bar{t}}
(y,\frac{\xi}{|\xi|}),p^{\bar{t}}(y,\frac{\xi}{|\xi|}))|\leq C\epsilon,\
j=1,\ldots,l \}.
%\label{oc3}
$$
That allows us to apply estimate \eqref{eq:rozcor} of
Lemma \ref{rozen:modify} with $q=0$ and $\delta=1$
to the operator $R_3$ to get
\be
\|R_3\|_{-3n/2-2}\leq C
\epsilon^{-n-1} \meas (\Xi),
\label{oc2}
\ee
where $\meas$ is the natural measure on $[0,T]^l\times S^*X$.
To use this inequality we need to estimate the measure
$\meas(\Xi)$. The set $\Xi$ maps to
$$
\Xi_2=\{(\bar{t},y,\xi):\ |T_j^0(\bar{t},
y,\xi)|\leq C\epsilon, \ |p^{\bar{t}}(y,\xi)|=1,\
j=1,\ldots,l\}
%\label{oc4}
$$
under the Hamiltonian flow $\Psi^{\bar{t}}$ which preserve
measure. The measure of $\Xi_2$ may be estimated by
the measure (in $[0,T]^l\times X\times\Rn$) of
$$
\Xi_3=\{(\bar{t},y,\xi):\ |T_j^0(\bar{t},
y,\xi)|\leq C\epsilon,\ c\leq |\xi|\leq C ,\ j=1,\ldots,l\}.
%\label{oc5}
$$
Let us now introduce sets
$$
\Sigma_j(t_1,\ldots,t_{j-1},y,\xi,\epsilon)=
\{t_j:\ |T_j^0(\bar{t},
y,\xi)|\leq C\epsilon,\ c\leq |\xi|\leq C \},\ j=1,\ldots,l.
%\label{oc6}
$$
To estimate measures of these sets we will use Theorem \ref{th:eps}.
For this, in the notation of Theorem \ref{th:eps},
we set $w=(y,\xi)$ and
consider the function
$$
t\mapsto f(t,w)=T_j^0(t_1,\ldots,t_{j-1},t,t_{j+1},\ldots,t_l,y,\xi),
$$
where $t$ takes the position of $t_j$ in $\bar{t}.$
Condition on Poisson brackets shows us that function
$f$ may have zeros in $t$ of order not greater then $M$.
It follows now from Theorem \ref{th:eps} that
\be
\max_{y,\xi}\meas\{ t_j\in [0,T]:\ T_j^0(\bar{t},y,\xi)\leq C\epsilon \}
\leq C\epsilon^{1/2M},
\label{eq:oc7}
\ee
where $\meas$ is just the Lebesgue measure, $y$ varies over a compact
set and $c\leq |\xi|\leq C$. Now we can use a simple observation that if
we have two functions $f_1(t_1)$ and $f_2(t_1,t_2)$ such that
$$
\meas\{t_1:\ |f_1(t_1)|\leq \epsilon\}\leq C\epsilon^{\alpha}, \
\max_{t_1}\meas\{t_2:\ |f_2(t_1,t_2)|\leq \epsilon\}\leq C\epsilon^{\alpha},
$$
then
\begin{eqnarray*}
\meas\{(t_1,t_2):\ |f_1(t_1)|\leq\epsilon, |f_2(t_1,t_2)|\leq\epsilon \}
& = &
\int_{|f_1(t_1)|\leq\epsilon, |f_2(t_1,t_2)|\leq\epsilon} dt_1 dt_2 \\
& = &
\int_{|f_1(t_1)|\leq\epsilon}\left(
\int_{|f_2(t_1,t_2)|\leq\epsilon} dt_2 \right) dt_1 \\
& \leq &
C\epsilon^{2\alpha}.
\end{eqnarray*}
Applying this argument $l$ times to the estimate \eqref{eq:oc7}, we get
$$
\meas(\Xi_3)\leq C\ \textrm{vol}\ (\supp_y b) \max_{(y,\xi)}
\prod_{j=1}^l \max_{t_1,\ldots,t_{j-1}}
\meas (\Sigma_j(\bar{t},y,\xi,\epsilon)),
$$
which implies
$$
\meas (\Xi)\leq C \epsilon^{l/2M}.
$$
Here we used that the support of the amplitude $b(\bar{t},x,y,\xi)$ is
compact with respect to $x$ and $y$, and that $T$ is small. 
Finally, combining this with estimate
\eqref{oc2}, we get
\be
\| R_3\|_{-3n/2-2}\leq C \epsilon^{l/2M-n-1}.
\label{estr3}
\ee
An application of Lemma \ref{l:interpolation} to
\eqref{r1}, \eqref{r2} and \eqref{estr3} yields the requred estimate
for $Q$ in Theorem \ref{th:q}.

\section{Various auxiliary results}

First we describe the lifting of the problem to insure that characteristic
roots of the principal symbol of $P$ are globally uniquely defined.
The proof follows the paper of Rozenblum \cite{Rozenblum}.
We give it here for the completeness and since we will need this construction
to determine the orders for $L^p$ estimates in Theorem \ref{th:lp}.
\begin{lemma}
Let $X$ be a smooth compact manifold without boundary of dimension $n\geq
3$. Then there exists a cover $\tilde{X}$ of $X$ with finitely many leaves,
such that on the lifting $\tilde{E}$ of $E$ to $T^*\tilde{X}$ branches 
of eigenvalues
$a_j(x,\xi)$ and eigenspaces $V_j(x,\xi)$ of the principal symbol $A(x,\xi)$
are smooth and globally well defined. The space
$L^2(\tilde{E})$ has a decomposition into a direct sum of $m$ spaces
such that
the matrix representation of $P(\tilde{x},\tilde{\xi})$
with respect to this decomposition
consists of pseudo-differential operators, and its principal symbol is
a diagonal matrix with $a_j(\tilde{x},\tilde{\xi})$ at the diagonal.
\label{l:diag}
\end{lemma}
\begin{proof}
Let us fix $z_0=(x_0,\xi_0)\in S^*X$. For each path in
$S^*X$ beginning at $z_0$ we look at the continuation of the
diagonalization along this path. At each point $z\in S^*X$ we
obtain up to permutation several possible collections $a_j(z), V_j(z)$,
such that homotopic paths from $z_0$ to $z$ correspond to the
same collection. Since for $n\geq 3$ homotopy of paths in $S^*X$
and their projections to $X$ is equivalent, we obtain a homomorphism
$\pi_1(X)$ to the permutation group of order $m$. A construction
in the homotopy theory (\cite{RF}) gives a cover
$p:\tilde{X}\to X$ with finitely many leaves, such that the lifting
of this homomorphism to $\pi_1(\tilde{X})$ is trivial. So, each
closed path in $S^*\tilde{X}$ takes eigenvalues and, therefore,
also eigenspaces, to themselves. This means that eigenvalues and
eigenspaces have global smooth branches on $S^*\tilde{X}$ and hence
also on $T^*\tilde{X}$.

Let $p_j^0(\tilde{x},\tilde{\xi})$ be a family of orthogonal
projectors on $V_j(\tilde{x},\tilde{\xi})$.
By the standard Gram-Schmidt process, we can add lower order terms
to $p_j^0(\tilde{x},\tilde{\xi})$ to obtain symbols
$p_j(\tilde{x},\tilde{\xi})$, for which
$p_j(\tilde{x},\tilde{D})p_k(\tilde{x},\tilde{D})=
\delta_{jk}p_j(\tilde{x},\tilde{D})$ and
$\sum p_j(\tilde{x},\tilde{D})=1$. For $u\in L^2(E)$, let
$u_j=p_j(\tilde{x},\tilde{D})u$.
Since $P(\tilde{x},\tilde{D}) p_j(\tilde{x},\tilde{D})=
p_j(\tilde{x},\tilde{D}) P(\tilde{x},\tilde{D})
p_j(\tilde{x},\tilde{D})$ modulo lower
order terms, $P=\sum p_j P p_k$ yields the desired diagonalization.
\end{proof}
Let $a_j$ still denote the smooth global branches of characteristics of
$A$ lifted to $\tilde{X}$. Let $\tilde{U}(t)$ be the fundamental
solution to \eqref{eq:eq1} with operator $P$ lifted to $\tilde{X}$. Let
$\tilde{U}(t,\tilde{x},\tilde{y})$ be the integral kernel of
$\tilde{U}(t)$, $\tilde{x}, \tilde{y}\in \tilde{X}$. Let
\be
U(t,x,y)=\sum_{p\tilde{y}=y} \tilde{U}(t,\tilde{x},\tilde{y}),
\ x,y\in X,
\label{eq:cover}
\ee
where $\tilde{x}$ is any point of $\tilde{X}$ such that $p\tilde{x}=x$
and the summation is carried out over all preimages $\tilde{y}$ of $y$.
Since \eqref{eq:cover} is invariant under permutations of the leaves
of the cover $\tilde{X}$, the kernel $U(t,x,y)$ is independent of the
choice of $\tilde{x}$. Equation and Cauchy data are automatically
satisfied, so \eqref{eq:cover} gives a fundamental solution to the
system \eqref{eq:eq1} on $X$.

The following Lemma gives some estimates for the operator norm from
$L^2$ to $H^s$ for Fourier integral operators in terms of $L^\infty$
norms of the amplitude and its derivatives. Recall that by
$\|T\|_s$ we denote the operator norm of $T$ from $L^2$ to $H^s$.
\begin{lemma}
Let $T$ be a Fourier integral operator
\be \label{eq:rozt}
Tu(x)=\int_\Rn\int_X e^{i\varphi(x,y,\xi)} a(x,y,\xi) u(y) dy d\xi,
\ee
where $X$ is an open set in $\Rn$. Assume that $a\in S^q$, $q\in\Z$,
is an amplitude of order $q$ and has compact support with respect to
$x$ and $y$ in $X$. Assume also that $\partial\varphi/\partial x\not =0$
for $\xi\not=0.$ Then $T$ extends to a bounded operator
from $L^2(X)$ to $H^{-q-n-1}(X)$ with
\be \label{eq:rozl}
||T||_{-q-n-1}\leq C||a\jp{\xi}^{-q}||_{C^{|q+n+1|}},
\ee
where $\jp{\xi}=(1+|\xi|^2)^{1/2}.$

Moreover, let us assume in addition that the support of $a$
belongs to a conical set with respect to $\xi$ which does not
depend on $x$, i.e.
$$
\supp a \subset \Xi_1=\{(x,y,\xi)
: x\in X, (y,\frac{\xi}{|\xi|})\in \Xi, \xi\not=0\},
$$
where $\Xi$ is subset of $S^*X$.
Then for any $\delta>0$ we have
\be \label{eq:rozcor}
||T||_{-q-3n/2-1-\delta}\leq C\meas(\Xi)||a\jp{\xi}^{-q}||_{C^{|q+n+1|}},
\ee
where $\meas$ is the canonical induced measure on $S^*X$ and constant $C$
may depend on the size of the support of $a$ with respect to $x$ and $y$.
\label{rozen:modify}
\end{lemma}

\begin{proof}
First we consider the case $q+n+1< 0.$ Differentiating
(\ref{eq:rozt}) $-(q+n+1)$ times with respect to $x$ we obtain
the integral with amplitude of order $-(n+1).$ The integral
with respect to $\xi$ converges absolutely, so
$T$ extends as a bounded operator from $L^2(X)$ to
$C^{-(q+n+1)}$ with estimate (\ref{eq:rozl}). Estimate
(\ref{eq:rozcor}) clearly follows from this as well.

We will now consider the case $q+n+1\geq 0.$
For a smooth function $v$, let us consider the bilinear form
$(Tu,v).$
Let us define operator $L$ as the transpose of
$^tL=(1+|\partial_x\varphi|^2)^{-1}
(1-i\partial_x\varphi\cdot\partial_x).$
Since $\partial\varphi/\partial x\not=0$,
integrating by parts $q+n+1$ times with operator $L$,
we obtain an absolutely
convergent integral with respect to $\xi$, and the estimate
$$ |(Tu,v)|\leq C\|u\|_{L^2} \|v\|_{q+n+1}
||a\jp{\xi}^{-q}||_{C^{q+n+1}}.$$
This implies (\ref{eq:rozl}).

We can slightly modify this argument to obtain an estimate of
operator $T$ acting from $L^\infty(X)$.
Indeed,
\begin{eqnarray*}
|(Tu,v)| & \leq & \left|\int_\Rn\int_\Rn\int_\Rn e^{i\varphi(x,y,\xi)}
L^{q+n+1}\left( a(x,y,\xi) \bar{v}(x)\right) u(y) dy d\xi dx\right| \\
& \leq & C\|v\|_{q+n+1}\|u\|_{L^\infty} I,
\end{eqnarray*}
where
$$
I^2=\int_X  I_1^2(x)dx, \;
I_1(x)=\int_X\int_{\mathbb{R}^n}|\tilde{a}(x,y,\xi)|d\xi dy,
$$
with some amplitude $\tilde{a}$ of order $-n-1$.
We can estimate $\tilde{a}$ by $\jp{\xi}^{-n-1}$ and
from the embedding theorems it follows that
$$ |(Tu,v)|\leq C \|v\|_{q+n+1}\|u\|_{n/2+\delta} \meas(\Xi)
\|a\jp{\xi}^{-q}\|_{C^{q+n+1}},$$
implying estimate (\ref{eq:rozcor}) for the norm of operator
$T$ acting on $L^2(X).$
\end{proof}

The following theorem shows that if a smooth function on
a bounded interval
has zeros only of finite order, then the measure of the set
where this function is small is also small. Moreover, if
we have a family of such functions continuously dependent on
a parameter, we can estimate measures of sets where functions
are small uniformly for all parameters varying over compact sets.
\begin{theorem}
Let $W\subset\Rn$ be compact and let $0<T<\infty.$
Let a real valued function $f=f(t,w)$ be continuous in $w\in W$
and smooth in $t\in [0,T]$ up to the boundary of $[0,T]$.
Let $M\in\N$ and
suppose that for each $w\in W$ function $f(\cdot,w)$
has zeros with respect to $t$ of order not greater than $M.$
Then there exist $C>0$ and $\epsilon_0>0$ such that for all
$\epsilon<\epsilon_0$ we have
$$ \sup_{w\in W} \meas\{t\in [0,T]: |f(t,w)|\leq\epsilon\} \leq
C\epsilon^{1/2M}.
$$
\label{th:eps}
\end{theorem}
\begin{proof}
For $\epsilon>0$ let
$$\Sigma(w,\epsilon)=\{t\in [0,T]:\; |f(t,w)|\leq C\epsilon\}.$$
Let $K(w)$
be the number of zeros of function $f(\cdot,w)$ with respect
to $t\in [0,T]$
and let $K$ be the maximum of $K(w)$ over $w\in W$.
It is obvious that $K$ is a finite number due to the condition
on zeroes of $f$ and compactness of $W$.

Let $\alpha_p>0$, $p\in\N$, be a decreasing sequence of positive
numbers which we will choose later. Let us define sets
$\Sigma^p(w,\epsilon)$ by setting
\ba
\label{sigm1}
\Sigma^p(w,\epsilon)=\{t\in [0,T]: |f(t,w)|\leq C\epsilon,\ldots, \\
|\partial_{t}^{p-1} f(t,w)|\leq
C\epsilon^{\alpha_{p-1}},\ |\partial_{t}^p f(t,w)|
\geq C\epsilon^{\alpha_p}\}, \ p\in\N. \nonumber
\ea
We claim now that there exists  $\epsilon_0>0$ such that for all
$0<\epsilon<\epsilon_0$, and all $w\in W,$
we have $\Sigma^p(w,\epsilon)=\emptyset$ for all $p>M.$
Indeed, if it is not so, then due to compactness of $W$ there are
sequences $t_n$, $w_n$ and $\epsilon_n$,
converging to $t^*\in [0,T]$, $w^*\in W$ and zero,
respectively,
such that
$$
|f(t_n, w_n)|\leq C\epsilon_n, \ldots ,
|\partial_{t}^{M+1} f(t_n, w_n)|\leq
C\epsilon_n^{\alpha_{M+1}},
$$
and consequently
$$
\partial_{t}^p f(t^*, w^*)=0, \ p=0,\ldots,M+1,
$$
which is impossible. It follows now that the set
$\Sigma(w,\epsilon)$ may be presented as the following
union of sets
\be
\Sigma(w,\epsilon)=\bigcup_{p=1}^M \Sigma^p(w,\epsilon).
\label{eq:decomp}
\ee
The idea of the proof now is to show first that the number of
connected components of sets $\Sigma^p(w,\epsilon)$ is finite
and can be estimated uniformly over all $w$ and $\epsilon.$
Then we will show that the size of each connected component is
small and can be estimated by $\epsilon^{1/2M},$ which will
imply Theorem \ref{th:eps}.
These statements are proved in the following two lemma.

\begin{lemma} There exists $\epsilon_0>0$
 such that for all
 $0<\epsilon<\epsilon_0$ the inequality
 \be
 \Delta(\Sigma^p(w,\epsilon))\leq K(M+1)^2,\; \forall w\in W,
 1\leq p\leq M,
 \label{oc10}
 \ee
 holds, where $\Delta(\Sigma^p)$ is the number of connected components of
 $\Sigma^p.$
 \label{l:components}
\end{lemma}
\begin{proof}
For simplicity let us first consider the case $p=1$.
Let us assume that \eqref{oc10} is not valid. Then there exist
sequences $w_n$ converging to some $w^*\in W$ and $\epsilon_n$ converging
to $0$ such that
\be
\Delta(\Sigma^1(w_n,\epsilon_n))> K(M+1)^2. \label{oc8}
\ee
 Let us now choose small enough $\epsilon_1>0$ such
 that each connected interval in the closure of
 $\Sigma(w^*,\epsilon_1)$ will
 include one and only one zero of function $f(t,w^*)$ with
 respect to $t$ and will not
 include zeros of derivative of $\partial_t f(t,w^*)$
 different from zeros of $f(t,w^*)$.
 This is possible because if we had an infinite number of
 zeros of $\partial_t f(t,w^*)$ approaching a zero $t^*$ of
 $f(t,w^*)$, it would mean that $\partial_t f(t^*,w^*)=0$ and
 that in fact $t^*$ is zero of $\partial^k f(t,w^*)$ for all
 $k\geq 1$, which is impossible since we assumed that all
 zeros of $f(t,w^*)$ are of finite order not exceeding $M$.
 It follows now that $\Delta(\Sigma(w^*,\epsilon_1))\leq K.$
 Since $w_n$ converges to $w^*$, by continuity we also have
 \be
  \Sigma^1(w_n,\epsilon_n)\subset\Sigma(w^*,\epsilon_1),
  \label{inter}
 \ee
 for sifficiently large $n$.
From \eqref {oc8} and \eqref{inter} it follows that there exist
$(M+1)^2+1\geq 2M+1$ connected components of $\Sigma^1(w_n,\epsilon_n)$ which
are all contained in one of the connected components of closure of
$\Sigma(w^*,\epsilon_1)$. Let us denote this connected component
of $\Sigma(w^*,\epsilon_1)$ by $I$.

According to the definition of sets $\Sigma^1(w_n,\epsilon_n)$,
between two of its connected components function $f(t,w_n)$
must become relatively large (i.e. $>C\epsilon$)
or its derivative
$\partial_t f(t,w_n)$ must become relatively small
(i.e. $<C\epsilon^{\alpha_1}$).
From this we see that
in the first case $\partial_t f(\cdot,w_n)$ must become zero
at some point between these components, while in the second case
$\partial^2_t f(\cdot,w_n)$ must become zero at some point there.
Since the number of components of $\Sigma^1(w_n,\epsilon_n)$
in $I$ is at least $2M+1$, it follows that at least one of these two
cases occurs at least $M+1$ times.

Let us consider these cases separately. In the first case,
for sufficiently large $n$, we have at least
$M$ zeros $s^1_n<\ldots <s^M_n$ of the derivative
$\partial_t f(t,w_n)$ contained in $I$.
It follows that $\partial^2_t f(\cdot,w_n)$ has at least $M-1$
different zeros in $I$, $\partial^3_t f(\cdot,w_n)$ has at least
$M-2$ different zeros in $I$, etc. In particular,
there are $\tau^k_n\in[s^1_n,s^M_n]$ such that
$\partial_t^k f(\tau^k_n,w_n)=0$ for $k=1,\ldots,M$.

Using compactness of $I$ and continuity of $f$, it follows that
there are subsequences $s^1_{n_i}, \ldots,s^M_{n_i}$, $\tau^k_{n_i}$
and $w_n$
converging to $s^1_*,\ldots,s^M_*$  and $\tau^k_*$, respectively, which
are all contained in $I$, and $w_n\to w^*$.
Since functions $f(t,w)$ are smooth
in $t$ we have
$$
\partial_t f(s^1_*,w^*)=0,\ldots,
\partial_t f(s^M_*,w^*)=0,\;  \textrm{and} \;
\partial_t^k f(\tau^k_*, w^*)=0, \; k=1,..,M.
$$
Since there are no zeros of derivative $\partial_t f(\cdot,w^*)$
on the interval $I$ except may be a point $t^*$ which is zero of function
$f(\cdot,w^*)$, we obtain $s^1_*=\ldots=s^M_*=t^*$.
From $\tau^k_n\in [s^1_n,s^M_n]$ it follows that
$\tau^k_*=t^*$ for all $k=1,\ldots,M$, which means
$$
 \partial_t^k f(t^*,w_*)=0\; \textrm{for all}\; k=0,\ldots,M.
$$
But this is impossible since function
$f(\cdot,w^*)$ may have zeros only of order $M$ in view of our assumption.

The second case is sightly different. Here,
for sufficiently large $n$, we have at least
$M$ zeros $s^1_n<\ldots <s^M_n$ of the second order derivative
$\partial^2_t f(t,w_n)$ contained in $I\backslash \Sigma^1(w_n,\epsilon_n)$.
Because in this case we assumed that $\Sigma^1(w_n,\epsilon_n)$ breaks
into at least $M+1$ components due to the failure of condition
$|\partial_t f(t,w)|\geq C\epsilon^{\alpha_1}$, it follows that
\be
 |\partial_t f(s^i_n,w_n)|\leq C\epsilon_n^{\alpha_1}, \ i=1,\ldots,M.
 \label{eq:s1z}
\ee
As $n$ tends to infinity, we can choose subsequences of $s^i_n$ converging
to some $s^i_*\in I$. Because $\partial_t f(\cdot,w^*)$ does not have
zeros in $I$ except may be for some $t^*\in I$ which is also the unique zero
of $f(\cdot,w^*)$ in $I$, we get that $s^1_*=\ldots=s^M_*=t^*.$ From
\eqref{eq:s1z} we also have $\partial_t f(t^*,w^*)=0.$

Now, to deal with higher order derivatives of $f$ at $t^*$,
similar to the first case, we get a collection of
$\tau^k_n\in[s^1_n,s^M_n]$ such that
$\partial_t^k f(\tau^k_n,w_n)=0$ for $k=2,\ldots,M+1$.
Because of compactness $\tau^k_n$ has a subsequence converging to some
$\tau^k_*$, $k=2,\ldots,M+1$. Again, we must have $\tau^k_*=t^*$
and hence also $\partial_t^k f(t^*,w^*)=0$,
for all $k=2,\ldots,M+1$. Since we already showed that
$f(t^*,w^*)=\partial_t f(t^*,w^*)=0$, this
contradicts the assumption that $f(\cdot,w^*)$
may have zeros of order up to $M$.

The argument for $p\geq 2$ is similar. For $t$ between two connected
components of $\Sigma^p$, at least one of conditions in \eqref{sigm1}
breaks down. Note that since we assumed that the total number of
components of $\Sigma^p(w_n,\epsilon_n)$ is larger than $K(M+1)^2$ and the
number of components of the larger set $\Sigma(w^*,\epsilon_1)$
is at most $K$,
we will have at least $(M+1)^2+1$ components of $\Sigma^p(w_n,\epsilon_n)$
in $I$. This means that these
$p+1$ conditions fail at least $(M+1)^2$ times.
Since $p+1\leq M+1$, there is a condition that will fail at least $M+1$ times.

In the case this is the last condition
$|\partial_t^p f(t,w)|\geq C\epsilon^{\alpha_p}$ that fails while conditions
$|\partial_t^i f(t,w)|\leq C\epsilon^{\alpha_i}$ remain valid, we can
argue similar to the second case of the argument with $p=1.$ In this
we need to have at least $M-p$ different zeros of
$\partial_t^{p+1} f(t,w)$, which is the case if we have at least
$M-p+1$ such components in $I$. This is true since we have at least
$M+1$ such components.

If one of the other conditions fails,
let us take the smallest $i$ for which condition
$|\partial_t^i f(t,w)|\leq C\epsilon^{\alpha_i}$ fails $M+1$ times.
Again, we need
to have at least $M-i$ different zeros of $\partial_t^{i+1} f(t,w)$,
which would follow from having $M-i+1$ such components in $I$. This
is again true since we have at least $M+1$ such components.
\end{proof}

Our next step is to show that each connected component of
$\Sigma^p(w,\epsilon)$, for $p=1,\ldots,M$, is small enough.

\begin{lemma}
The length of each connected component of $\Sigma^p(w,\epsilon)$, for
$p=1,\ldots,M$, is no greater than $C\epsilon^{\alpha_{p-1}-\alpha_{p}}$.
\label{l:size}
\end{lemma}
\begin{proof}
Let $I$ be a connected component of $\Sigma^p(w,\epsilon)$
and let  $t^*\in I$. We are going to
estimate the maximal shift $\delta_0>0$ such that $t^*+\delta \in I$
for all $0<\delta<\delta_0$. Recalling the definition of
$\Sigma^p(w,\epsilon)$, we see that
\be
|\partial_t^{p-1} f(t^*+\delta,w)|\leq C\epsilon^{\alpha_{p-1}},\
|\partial_t^p f(t^{*}+\delta,w)|\geq
C\epsilon^{\alpha_p}.
\label{eq:estt1}
\ee
Using the Taylor expansion of $\partial_t^{p-1} f(\cdot,w)$ at
$t^*$, we have
$$
\partial_t^{p-1} f(t^*+\delta,w) = \partial_t^{p-1} f(t^*,w) +
\partial_t^{p} f(t^{**},w)\delta,
$$
where $t^{**}$ is some point between $t^*$ and $t^*+\delta.$
Since $t^*, t^*+\delta\in I$, it follows that $t^{**}\in I$ and
$$
|\partial_t^{p} f(t^{**},w)| |\delta|\leq
|\partial_t^{p-1} f(t^*+\delta,w)| + |\partial_t^{p-1} f(t^*,w)|
\leq 2C\epsilon^{\alpha_{p-1}}.
$$
From this and \eqref{eq:estt1} we obtain that
$|\delta_0|\leq C\epsilon^{\alpha_{p-1}-\alpha_{p}}$.
Consequently, the length of each connected component of
$\Sigma^p(w,\epsilon)$ can be estimated by
$C\epsilon^{\alpha_{p-1}-\alpha_{p}}$.
\end{proof}

Now we can finish the proof of Theorem \ref{th:eps}.
Let us choose $\alpha_k=1-k/2M$, $k=0,\ldots,M$.
According to Lemma \ref{l:size}
the length of each connected component of $\Sigma^p(w,\epsilon)$
can be estimated by $C\epsilon^{1/2M}$. Then according to Lemma
\ref{l:components}, the size of $\Sigma^p(w,\epsilon)$ can be
estimated by $CK(M+1)^2\epsilon^{1/2M}$. Because of
decomposition \eqref{eq:decomp} the size of $\Sigma(w,\epsilon)$
is estimated by $C\epsilon^{1/2M}$. Statement of Theorem
\ref{th:eps} is now a consequence of continuity of $f$ with respect
to $w$ and compactness of $W$.
\end{proof}

The following interpolation lemma shows that if
a function $u$ can be decomposed for all sufficiently
small $\epsilon$
into a sum of $u_1^{(\epsilon)}+u_2^{(\epsilon)}$ with a good
estimate for the norm of $u_2^{(\epsilon)}$ in a ``bad'' Sobolev
space $H^p$ with small index $p$, and with a bad estimate for
the norm of $u_1^{(\epsilon)}$ in a ``good'' Sobolev space
$H^r$ with large index $r$, then $u$ belongs to some intermediate
space $H^q$ with $p<q<r.$

\begin{lemma} \label{l:interpolation}
Let $u\in H^p(\Rn)$. Suppose that for all small enough
$\epsilon$ there is a
representation $u=u_1^{(\epsilon)} +u_2^{(\epsilon)}$ such that
\be
\|u_1^{(\epsilon)} :H^r\|\leq C\epsilon^{-T},\ \
\|u_2^{(\epsilon)} :H^p\|\leq C\epsilon^{S}, \; p<r; \; S, T>0.
\label{oc44}
\ee
Then $u\in H^q(\Rn)$ for any
\be
q<(pT+rS)(T+S)^{-1}.
\label{oc46}
\ee
\end{lemma}
\begin{proof}
Let a sequence $\{a_j\}$ be such that $0=a_0<a_1<\ldots$ and
$a_j\rightarrow +\infty$ as  $j\rightarrow +\infty$, and let
a sequence $\{b_j\}$ of positive numbers $b_j>0$ tend to zero.
Let us assume that $u\in H^q(\Rn)$ for some $q$. Then
\begin{eqnarray*}
\|u:H^q\|^2 & = & \int_{\Rn}|\hat{u}|^2 (1+|\xi|^2)^{q} d\xi \\
 & = & 2\sum_{j=0}^{+\infty}
\left(\int_{a_j\leq |\xi|\leq a_{j+1}} |\widehat{u^{(b_j)}_2}|^2
(1+|\xi|^2)^{q}d\xi
+\int_{a_j\leq |\xi|\leq a_{j+1}} |\widehat{u^{(b_j)}_1}|^2
(1+|\xi|^2)^{q} d\xi\right) \\
 & \leq & 2\sum_{j=0}^{+\infty} \left(
 \max \{\jp{a_j}^{2(q-p)},\jp{a_{j+1}}^{2(q-p)}\}
\int_{a_j\leq |\xi|\leq a_{j+1}}
|\widehat{u^{(b_j)}_2}|^2 (1+|\xi|^2)^{p} d\xi \right. \\
 & + & \left. \max\{\jp{a_j}^{2(q-r)}, \jp{a_{j+1}}^{2(q-r)}\}
\int_{a_j\leq |\xi|\leq a_{j+1}}
|\widehat{u^{(b_j)}_1}|^2(1+|\xi|^2)^{r}d\xi \right) .
\end{eqnarray*}
It follows from this estimate and \eqref{oc44} that
\be
\|u:H^q\|^2 \leq
C\sum_{j=0}^{+\infty}\left(
\max\{\jp{a_j}^{2(q-p)},\jp{a_{j+1}}^{2(q-p)}\}b_j^{2S}+
\max\{\jp{a_j}^{2(q-r)},\jp{a_{j+1}}^{2(q-r)}\}b_j^{-2T}
\right).
\label {oc45}
\ee
Now we are going to demonstrate that under
hypothesis \eqref{oc46} we may choose sequences  $\{a_j\}$ and
$\{b_j\}$ in such way that the right hand side of \eqref{oc45} will
be finite, so that conclusion of Lemma \ref{l:interpolation}
will follow. Let us set $a_j=j^{\alpha}$ and
$b_j=j^{-\beta}$ with some $\alpha, \beta>0$.
Then the series in \eqref{oc45} will
converge if the following inequalities are fulfilled
$$
2\alpha(q-p)-2\beta S<-1, \ \ 2\alpha(q-r)+2\beta T<-1.
$$
These inequalities can be transformed into
$$
\alpha(2T(q-p)+2S(q-r) )<-T-S,\ \ 2\beta N>
\alpha2(q-p)+1,
$$
which hold with positive constants $\alpha$ and $\beta$ if and only
if \eqref{oc46} is valid.
\end{proof}

\section{Spectral asymptotics}

In this section we will prove Theorem \ref{th:spectral}. Let $X$ be
a smooth compact manifold without boundary of dimension $n\geq 3$.
Then an elliptic operator $P(x,D)$ has a collection of eigenfunctions
and eigenvalues $\lambda_j\to\infty$. We will be interested in
distribution of eigenvalues and will find the asymptotics for
the spectral function $N(\lambda)=\sharp\{j:\lambda_j<\lambda\}$.
One of the most effective methods to study such asymptotics is
to use an explicit representation for the fundamental solution of
the corresponding hyperbolic problem \eqref{eq:eq1}. We will follow
the method developed in \cite{Hormander-sp}, \cite{DG},
\cite{SV}, etc., for scalar operators, and in \cite{Rozenblum}
for the case $M=1$ in Condition C.
The following proposition was implicitely proved
in \cite{DG} and formulated in \cite{Rozenblum}.

\begin{prop}
Let $\chi_1, \chi_2\in C_0^\infty(\R)$ be such that
$\chi_1(0)=1$ and $0\not\in\supp\chi_2$. Suppose that
\be
\begin{gathered}
\tr \Fcal^{-1}_{t\to\mu}(\chi_1(t)U(t))=c_1\mu^{n-1}+c_2\mu^{n-2}
+o(\mu^{n-2}), \\
\tr \Fcal^{-1}_{t\to\mu}(\chi_2(t)U(t))=o(\mu^{n-1}), \ \mu\to\infty.
\end{gathered}
\label{eq:f1}
\ee
Then $N(\lambda)=c_1 n^{-1}\lambda^n+c_2(n-1)^{-1}\lambda^{n-1}
+o(\lambda^{n-1}).$
\label{prop:as}
\end{prop}
In previous sections, and in particular in the proof of Theorem
\ref{th:lp}, we have represented the propagator $U(t)$ in \eqref{eq:eq1}
as an infinite sum of some extensions of Fourier integral operators, i.e.
\be
U(t)=\sum_{j=0}^\infty e^{-i\tilde{A}t} Q_j,
\label{eq:rep1}
\ee
where $Q_0=I$ and $\tr Q_1=0$. So we also have 
\be
\tr e^{-i\tilde{A}t} Q_1=0.
\label{eq:tr1}
\ee
Let us first consider the contribution of the first term, which is
$e^{-i\tilde{A}t}$. It is the propagator for a block-diagonal system,
so the asymptotic behaviour of 
$\Fcal^{-1}_{t\to\mu}(\chi_\sigma(t)e^{-i\tilde{A}t})$, $\sigma=1,2$, 
determines the
spectral distribution for a system of independent scalar equations.
Thus, this is the sum of spectral distributions for scalar operators,
which are well-known (e.g. \cite{Hormander-sp}, \cite{DG}). So
$$
\Fcal^{-1}_{t\to\mu}\tr(\chi_1(t)e^{-i\tilde{A}t})=
c_1\mu^{n-1}+c_2\mu^{n-2}+o(\mu^{n-2}).
$$
We will show that under suitable conditions the asymptotics in
\eqref{eq:f1} are determined only by the first term of the
series \eqref{eq:rep1}. Let us first observe that in view of
Theorem \ref{th:q} operator 
$\sum_{j=K}^\infty e^{-i\tilde{A}t} Q_j$ is compact for sufficiently
large $K$, so it does not contribute to the singularity of
$\tr \chi_\sigma U$, $\sigma=1,2$. Therefore, we can replace $U(t)$ in
\eqref{eq:f1} by a sum of finitely many terms in \eqref{eq:rep1}.
Also, in view of \eqref{eq:tr1}, the second term of
\eqref{eq:rep1} does not contribute to \eqref{eq:f1}. Other
terms of the sum \eqref{eq:rep1} are of the form
$$
e^{-i\tilde{A}t}\int_0^t\int_0^{t_1}\cdots\int_0^{t_l}
Z(t_1)Z(t_2)\cdots Z(t_{l+1}) dt_{l+1}\ldots dt_1,\ l\geq 2,
$$
where $Z(\tau)=-ie^{i\tilde{A}\tau}B e^{-i\tilde{A}\tau}$.
In this way, after a change of variables $s_1=t-t_1$, 
$s_2=t_1-t_2, \ldots, s_{l+1}=t_l$, 
we obtain a sum of terms of the form
$$
I_\sigma(\mu)=\tr \int \chi_\sigma(s_1+\cdots+s_{l+1}) 
e^{i(s_1+\cdots+s_{l+1})\mu}
B(s)L(s)ds, \ \sigma=1,2.
$$
where $B(s)\in\Psi^0$,
$L(s)=e^{-i\tilde{A}_{j_1}s_1} e^{-i\tilde{A}_{j_2}s_2}\ldots
e^{-i\tilde{A}_{j_{l+1}}s_{l+1}}$, and $\tilde{A}_j$ is the $j$-th block
of $\tilde{A}$ corresponding to $a_j$. As before, $L(s)$ may be represented
as a locally finite sum of oscillatory integrals with phase
$\varphi(s,x,y,\xi)$, for which
$$
\frac{\partial\varphi}{\partial s_k}=S_{j_k}(s,x,\frac{\partial\varphi}
{\partial x}),
$$
where $S_{j_k}(s,x,p)=a_{j_k}(\Phi_{j_{k-1}}^{s_{k-1}}\cdots
\Phi_{j_1}^{s_1}(x,p))$. Therefore,
$$
I_\sigma(\mu)=\int \int_X\int_\Rn \chi_\sigma(s_1+\cdots+s_{l+1})
e^{i(s_1+\cdots+s_{l+1})\mu+i\varphi(s,x,x,\xi)}b(s,x,x,\xi) d\xi dx ds.
$$
Substituting 
$
\xi=\mu\tau\omega,\ \tau>0, \ |\omega|=1,
$
we get
\be
I_\sigma(\mu)=\mu^n\int\int_X \int_{{\bf S}^{n-1}}\int_0^\infty
\chi_\sigma(s_1+\cdots+s_{l+1})e^{i\mu(s_1+\cdots+s_{l+1}+\tau\varphi)}
b\tau^{n-1} d\tau d\omega dx ds.
\label{eq:im1}
\ee
To finish the argument for this operator, we can follow 
\cite{Rozenblum} to show that smoothing of $Q_l$'s implies
\eqref{eq:f1}. If we change variables again by
$\rho=\sum s_j, s_j=\rho\kappa_j$ and introduce
$K=\{\kappa_j\geq 0, \ \sum\kappa_j=1\}$, we get
$$
I_1(\mu)=\mu^n\int_K\int_X\int_{{\bf S}^{n-1}}\int_{-\infty}^\infty
\int_0^\infty \chi_1(\rho) e^{i\mu(\rho+\tau\varphi(\rho\kappa,x,x,\omega)}
\rho^l\tau^{n-1} b d\tau d\rho d\omega dx d\kappa.
$$
For fixed $\omega,x,\kappa$ the point 
$\tau=-(\sum\kappa_k S_{j_k}(\rho\kappa,x,\partial\varphi/\partial x))^{-1}$,
$\rho=1$, is a nondegenerate stationary point of this oscillatory integral,
implying $I_1(\mu)=O(\mu^{n-1-l})=o(\mu^{n-2})$. This already gives H\"ormander's
first term of $N(\lambda)$. Now we will
show that singularities at $t\not=0$ do not
give essential contributions to the second term of
spectral asymptotics. In the
analysis of $\tr \chi_2 U$, the contribution of the
first term of \eqref{eq:rep1} was established in \cite{DG}. 
Let us look at the integral \eqref{eq:im1}
with respect to $s_1, s_2,\tau$ with
fixed $x,\omega,s_3,\ldots,s_{l+1}$. The stationary phase method with
repsect to $\tau, s_1$ gives a stationary point
$\varphi(s,x,x,\omega)=0$, $\tau=-(\partial\varphi/\partial s_1)^{-1}$.
It is non-degenerate because 
$\det \partial_\tau\partial_{s_1}(\tau\varphi)=-a_1(s,x,\partial_x\varphi)^2$.
This gives the estimate $I_2(\mu)=O(\mu^{n-1})$. Since 
$\partial\varphi/\partial s_2=0$ only at 
$\tau=-(\partial\varphi/\partial s_2)^{-1}=a_2(s,x,\partial_x\varphi)^{-1}$,
the phase is stationary with respect to $s_2$ only at $(x,\xi)$ for which
$a_1(s,x,\partial_x\varphi)=a_2(s,x,\partial_x\varphi)$. This is the set
of measure zero and, therefore, $I_2(\mu)=o(\mu^{n-1})$, which shows
\eqref{eq:f1}. Finally we note that if the support of $\chi_1$ is sufficiently
small (i.e. when $s$ is small), 
terms in \eqref{eq:cover}, corresponding to $\tilde{x}$ and $\tilde{y}$
in different leaves of $\tilde{X}$ will produce only smoothing 
operators \eqref{eq:cover}
in view of the finite propagation speed. So their contribution
to $I(\mu)$ is rapidly decreasing and we obtain 
the result also on $X$.

%%%%%%%%%%%%%%%%%%%%%%%%%%%%%references%%%%%%%%%%%%%%%%%%%%%%%%%%%%%%%%%%%
%\newpage

\end{document}